\documentclass[11pt]{article}
\usepackage[T1]{fontenc}
\usepackage[utf8]{inputenc}
\usepackage{amsfonts}
\usepackage{amsmath}
\usepackage{enumerate}
\usepackage{amssymb}
\usepackage{amsthm}
\usepackage{bbm}
\usepackage{bm}
\usepackage{mathrsfs}
\usepackage{verbatim}
\usepackage{setspace}
\usepackage{color}
\usepackage{pdfsync}
\usepackage{booktabs}
\usepackage{tikz-cd}

\usepackage{tikz}

\usepackage{graphicx}
\usepackage{amsmath,varwidth,array,ragged2e}
\usepackage{tabu}
\usepackage{array}
\usepackage{makecell}

\usepackage{enumitem}
\theoremstyle{plain}
\newtheorem{theorem}{Theorem}[section]
\newtheorem{proposition}[theorem]{Proposition}
\newtheorem{lemma}[theorem]{Lemma}

\theoremstyle{definition}

\newtheorem{remark}[theorem]{Remark}

\newtheorem{example}[theorem]{Example}

\newtheorem{assumption}[theorem]{Assumption}

\theoremstyle{remark}

\renewenvironment{thebibliography}[1]{%
\begin{oldthebibliography}{#1}%
\setlength{\baselineskip}{.9em}
\linespread{1}
\small
\setlength{\parskip}{0.3ex}%
\setlength{\itemsep}{.5em}%
}%
{%
\end{oldthebibliography}%
}

\DeclareMathOperator{\id}{Id}

\DeclareMathOperator*{\argmin}{arg\, min}

\numberwithin{equation}{section}

\usepackage[pdfborder={0 0 0}]{hyperref}
\hypersetup{
  urlcolor = black,
  pdfauthor = {Stephan Eckstein, Michael Kupper},
  pdfkeywords = {T.B.D.;
},
}

\usepackage{arydshln}

\begin{document}

\title{\vspace{-0em}
	Computation of optimal transport and related hedging problems via 
	penalization 
	and neural networks\thanks{We thank Daniel Bartl, Fabian Both, Jens 
	Jackwerth, 
		Mathias Pohl, Stefan Volkwein and two anonymous referees for helpful comments.}
	\date{\today}
	\author{
		Stephan Eckstein%
		\thanks{
			Department Mathematics and Statistics, University of Konstanz, 
			\texttt{stephan.eckstein@uni-konstanz.de}.
		}
		\and
		Michael Kupper%
		\thanks{
			Department Mathematics and Statistics, University of Konstanz, 
			\texttt{kupper@uni-konstanz.de}.
		}
	}
}
\maketitle \vspace{-1.2em}

\begin{abstract}
	This paper presents a widely applicable approach to solving 
	(multi-marginal, martingale) optimal transport and related problems via 
	neural networks. The core idea is to penalize the optimization problem 
	in its dual formulation and reduce it to a finite dimensional one which 
	corresponds to optimizing a neural network with smooth objective function. 
	We present numerical examples from optimal transport, martingale 
	optimal transport, portfolio optimization under uncertainty and generative 
	adversarial networks that showcase the generality and effectiveness of the 
	approach.
\end{abstract}

\vspace{.9em}

{\small
	
	
	\section{Introduction}
	
	In this paper we present a penalization method which allows to compute a 
	wide class of optimization problems
	of the form
	\[
	\phi(f) = \sup_{\nu\in \mathcal{Q}} \int f\,d\nu
	\]
	by means of neural networks. The most widely known representative of such a functional occurs in the optimal transport problem, to be introduced shortly.
	More generally, these functionals appear for instance in the representation
	of coherent risk measures \cite{artzner1999coherent} as the worst-case expected 
	loss over a class $\mathcal{Q}$ of scenario probabilities, in the 
	representation of nonlinear expectations \cite{peng2007g}, or as the upper 
	bound of arbitrage-free prices for a contingent claim $f$, see 
	e.g.~\cite{follmer2011stochastic}. To solve the initial problem $\phi(f)$ 
	we will make use of its dual formulation and restrict to the subclass of 
	those optimization problems which can be realized as a minimal superhedging 
	price
	\[\phi(f)=\inf_{\substack{h \in 
			\mathcal{H}:\\ h \geq f}} \int h \,d\mu_0.\] 
	for some $\mu_0 \in \mathcal{Q}$, where $\mathcal{H}$ is a set of continuous and bounded functions $h: \mathcal{X} \rightarrow \mathbb{R}$, where the relation of $\mathcal{H}$ and $\mathcal{Q}$ is given at the beginning of Section \ref{sec2}.
	A very similar class of
	optimization problems in an abstract framework of Banach lattices is 
	studied in \cite{ekren2017constrained}. 
	Under sufficient regularity conditions the values of the primal problem 
	$\sup_{\nu\in \mathcal{Q}} \int f\,d\nu$ and its dual problem 
	$\inf_{\substack{h \in \mathcal{H}:\,h \geq f}} \int h \,d\mu_0$ can be shown to coincide, see 
	e.g.~\cite{cheridito2017duality} for related pricing-hedging dualities.  
	
	A typical example is the Kantorovich relaxation 
	\cite{kantorovich1942translocation}
	of Monge's optimal transport problem, where $\mathcal{Q}$ is the set of 
	probability measures on a product space $\mathcal{X} = 
	\mathcal{X}_1\times\mathcal{X}_2$ with given marginals $\mu_1$ and $\mu_2$, 
	and where $\mathcal{H}$ is the set of all continuous and bounded functions 
	$h(x_1, x_2) = h_1(x_1) + h_2(x_2)$ and $\int h \,d\mu_0 = \int_{\mathcal{X}_1} h_1 
	\,d\mu_1 + \int_{\mathcal{X}_2} h_2 \,d\mu_2$. Further frequently studied 
	problems in this class
	include multi-marginal optimal transport and Wasserstein distances (see 
	e.g.~\cite{bartl2017duality, vallender1974calculation, 
	villani2008optimal}), 
	martingale optimal transport (see e.g.~\cite{beiglbock2013model, 
		galichon2014stochastic, guo2017computational, henry2013automated}), 
		value at 
	risk under dependence uncertainty (see e.g.~\cite{bernard2017robust, 
		embrechts2013model, puccetti2012computation}), or calculating worst 
		case copula 
	values and improved Fr\'echet-Hoeffding bounds (see 
	e.g.~\cite{bartl2017sharpness, lux2017improved}). Moreover, $\phi(f)$ 
	serves as 
	a building block for several other problems, like generative adversarial 
	networks (where additionally, the 
	optimization includes generating a distribution, see 
	e.g.~\cite{arjovsky2017wasserstein,feizi2017understanding, 
		gulrajani2017improved}), portfolio choice under 
	dependence uncertainty (where 
	additionally, portfolio weights are optimized, see 
	e.g.~\cite{bernard2014risk, 
	pflug2017review}), or robust optimized certainty equivalents (see e.g.~\cite{eckstein2018robust}). In these cases, the solution approach presented 
	in 
	this paper is still applicable.

	\paragraph{Summary of the approach} The goal is to solve $\phi(f)$ 
	numerically. 
	The implementation will build on the dual representation of $\phi(f)$. 
	The first step is to go over to a finite dimensional setting, where the set $\mathcal{H}$ is replaced by a subset $\mathcal{H}^m$:
	\[
	\phi^m(f) = \inf_{\substack{h\in\mathcal{H}^m:\\ h\geq f}} \int h \,d\mu_0
	\]
	Theoretically, we will look at a sequence 
	$(\mathcal{H}^m)_{m\in\mathbb{N}}$ with $\mathcal{H}^1 \subseteq 
	\mathcal{H}^2 \subseteq ... \subseteq \mathcal{H}$ such that 
	$\mathcal{H}^{\infty} := \cup_{m \in \mathbb{N}} \mathcal{H}^m$ is in a 
	certain sense dense in $\mathcal{H}$.
	More concretely, $\mathcal{H}^m$ can be a set of neural networks with a fixed structure (but unspecified parameter values), and $m$ measures the number of neurons per layer.
	
	To allow for a step-wise updating of the parameters (e.g.~by gradient descent methods) for the space $\mathcal{H}^m$, the inequality constraint $h \geq f$ is penalized. To this end, we introduce a \textsl{reference probability measure} $\theta$ on the state space $\mathcal{X}$. Intuitively, this measure will be used to sample points at which the inequality constraint $h\geq f$ can be tested. Further, we introduce a differentiable and nondecreasing \textsl{penalty function} $\beta : \mathbb{R} \rightarrow \mathbb{R}_+$. This leads to the penalized problem
	\begin{align*}
		\phi_{\theta, \beta}^m(f) &= \inf_{h\in\mathcal{H}^m} \Big\{ \int h \,d\mu_0 + \int \beta(f-h) \,d\theta \Big\}.
	\intertext{For theoretical considerations we also introduce}
		\phi_{\theta, \beta}(f) &= \inf_{h\in\mathcal{H}} \Big\{ \int h \,d\mu_0 + \int \beta(f-h) \,d\theta \Big\}.
	\end{align*}
	Theoretically, we will again consider sequences of penalty functions $(\beta_{\gamma})_{\gamma > 0}$ parametrized by a penalty factor $\gamma$, and use the notation $\phi_{\theta, \gamma}(f) := \phi_{\theta, \beta_\gamma}(f)$ and $\phi^m_{\theta, \gamma}(f) := \phi^m_{\theta, \beta_\gamma}(f)$. Here, an increasing penalty factor can be seen as a more and more precise enforcing of the inequality constraint $h\geq f$.
	
	The problems $\phi^m_{\theta, \gamma}(f)$ are the ones which are solved numerically. Chapters \ref{sec2} and \ref{sec3} study the relation between this problem which is eventually implemented, and the initial problem $\phi(f)$. To this end, we analyse how the introduced approximative problems behave for $m \rightarrow \infty$ and $\gamma \rightarrow \infty$. Figure \ref{figureoverview} summarizes the occurring problems and their relations. Notably, we are only interested in convergence of optimal values, not that of optimizers.
	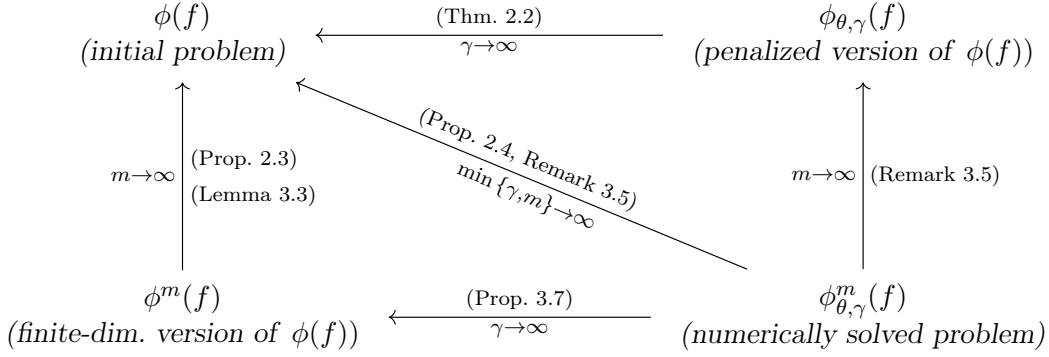
\begin{figure}
		\label{figureoverview}
		\begin{equation*}
		\begin{tikzcd}[row sep=2.5cm, column sep=3.5cm]
		\begin{tabular}{c} $\phi(f)$ \\ \textsl{(initial problem)} \end{tabular} & \begin{tabular}{c} $\phi_{\theta, \gamma}(f)$ \\ \textsl{(penalized version of }$\phi(f)$) \end{tabular} \arrow[l, "\gamma \rightarrow \infty", "(\text{Thm.~}\ref{thm:reg})"{above}] \\
		\begin{tabular}{c} $\phi^m(f)$ \\ \textsl{(finite-dim.~version of }$\phi(f))$ \end{tabular} \arrow[u, "m \rightarrow \infty", "{\hspace{-2mm}\begin{tabular}{l} \scriptsize{(Prop.~\ref{prop:approx})} \\ \scriptsize{(Lemma \ref{lem:CondD})} \end{tabular}}"{right}] &  \begin{tabular}{c} $\phi^m_{\theta, \gamma}(f)$ \\ \textsl{(numerically solved problem)} \end{tabular} \arrow[u, "m \rightarrow \infty", "(\text{Remark \ref{uniform}})"{right}] \arrow[l, "\gamma \rightarrow \infty", "(\text{Prop.~}\ref{nouniform})"{above}] \ar[ul, "\min{\{\gamma, m\}} \rightarrow \infty~~~"{sloped, near start}, "(\text{Prop.~\ref{prop:unif}, Remark \ref{uniform}})"{sloped, above}, anchor=north west]
		\end{tikzcd}
		\end{equation*}
		\caption{Occurring problems and their relations. The depicted convergences are studied in Section \ref{sec2} and, in a more specific context of neural networks, in Section \ref{sec3}.}
	\end{figure}
	
	The final step is to find a numerical solution of $\phi_{\theta, 
		\gamma}^m(f)$, 
	which means in practice finding the optimal parameters of the network 
	$\mathcal{H}^m$. We use Tensorflow \cite{tensorflow2015-whitepaper} and the 
	Adam optimizer \cite{kingma2014adam} to this end, and thus mostly regard 
	this step as a black box. We will denote the numerical optimal solution by 
	$\hat{\phi}^m_{\theta,\gamma}(f)$.

	\paragraph{Implementation method: Related literature}  Penalization of optimal transport problems 
	has 
	been studied in several works (see e.g.~\cite{benamou2015iterative, 
		carlier2017convergence, cominetti1994asymptotic, cuturi2013sinkhorn, 
		NIPS2016_6566, gulrajani2017improved, schmitzer2016stabilized, 
		seguy2017large, solomon2015convolutional}). Entropic penalization in 
		particular 
	is applied often, which is in close relation to the Schrödinger problem 
	\cite{leonard2013surveypub}. Cominetti and San Mart{\'\i}n's work 
	\cite{cominetti1994asymptotic} from 1994 on entropic penalization of 
	arbitrary 
	linear programs can be applied to purely discrete optimal transport. The 
	basic 
	idea in \cite{cominetti1994asymptotic} is to obtain a strictly convex 
	problem 
	through penalization which can be solved quicker and converges to the 
	initial 
	problem, for an increasing penalty factor. More recently, Cuturi 
	\cite{cuturi2013sinkhorn} gives an efficient algorithm to compute discrete 
	optimal transport problems with two marginals based on entropic 
	penalization 
	and Sinkhorn's matrix scaling algorithm. Genevay et 
	al.~\cite{NIPS2016_6566} and Solomon et 
	al.~\cite{solomon2015convolutional} go further in this direction and give 
	algorithms to compute arbitrary optimal transport problems with two 
	marginals, 
	where the algorithm (for the case of continuous marginals) is based on a 
	reproducing kernel Hilbert space approach, and discretization, 
	respectively.   
	In \cite{NIPS2016_6566} the authors already mention that more general 
	regularizations beyond the entropic one are possible.
	Among others Benamou et al.~\cite{benamou2015iterative} and Schmitzer 
	\cite{schmitzer2016stabilized} use scaling algorithms related to 
	\cite{cuturi2013sinkhorn} for a larger class of problems, including for 
	example (discrete) multi-marginal, constrained and unbalanced optimal 
	transport. Carlier et al.~\cite{carlier2017convergence} show $\Gamma$-convergence of the entropic penalized 
	Wasserstein-2 distance to the unpenalized one. The same kind of $\Gamma$-convergence is also subject of the studies related to the Schrödinger problem \cite{leonard2013surveypub}, even for more general cost functions.  Recent research by Arjovsky 
	et al.~\cite{arjovsky2017wasserstein, gulrajani2017improved} inspired by 
	generative adversarial networks include solving a particular optimal 
	transport problem (the Wasserstein-1 distance) based on $L^2$ penalization. In these works, the numerical approach to solve optimal transport problems by parametrization of the dual variables by neural networks originated. Seguy et al.~\cite{seguy2017large} apply a neural 
	network based approach to arbitrary optimal transport problems with two 
	marginals. Their theoretical results are broadly based on entropic 
	penalization, discretization, and weakly continuous dependence of the 
	optimal transport problem on the marginals.
	
	\paragraph{Contribution} The current paper gives a unifying numerical 
	solution 
	approach to problems of the form $\phi(f)$ based on penalization and neural 
	networks. The focus lies both on general applicability with respect to the 
	choice of problem, and also on a flexible framework regarding the solution 
	method.
	
	Compared to the existing literature, which often focusses on a single representative (often the optimal transport problem) among problems of the form $\phi(f)$, our theoretical results are widely applicable. Similarly, the penalization method and the resulting dual relations in this paper allow for many different forms of reference measure $\theta$ and penalty function $\beta_\gamma$, while the existing literature is often restricted to uniform or product reference measures, and exponential penalty functions.\footnote{In discrete settings, the 
		reference measure is usually the uniform distribution (see 
		e.g.~\cite{cuturi2013sinkhorn}, where the penalization does not 
		explicitly 
		include a reference measure. The penalization applied is simply the 
		entropy of 
		a measure, with corresponds to the relative entropy with uniform 
		reference 
		measure). In non-discrete settings, usually the product measure of the marginals specified by the 
		optimal 
		transport problems are used (see e.g.~\cite{NIPS2016_6566, seguy2017large}).} 
	We show the effects of 
	different reference measures and different penalty functions both 
	theoretically in Theorem \ref{thm:reg} and practically in 
	the numerical examples in Section \ref{sec:NumEx}. In some examples the 
	choice of an appropriate reference measure is crucial, see e.g.~Section 
	\ref{subsec::PO}.
	Equation \eqref{dualoptimizer} of Theorem \ref{thm:reg} also motivates an updating procedure for the reference measure to reduce the error arising from penalization, which is applied in Section \ref{sec::Rearrangement}.
	
	The presented approach is showcased with several examples, which are 
	mostly toy problems taken from existing papers. The reason we use toy 
	problems 
	is to allow for an evaluation of the numerical methods that can be based on analytical 
	solutions.
	
	\paragraph{Structure of the paper} In Section \ref{sec2} we present the 
	theoretical results on approximation
	and regularization. Section \ref{sec3} discusses the particular case of 
	$\mathcal{H}^m$ as built by multilayer feedforward networks. In Section 
	\ref{sec:NumEx} we illustrate the proposed method with several examples. 
	All 
	proofs are postponed to Section \ref{sec5}.

\section{Regularization and approximation of hedging functionals}\label{sec2}
Let $\mathcal{P}(\mathcal{X})$ be the set of all Borel 
probability measures on a Polish space $\mathcal{X}$, and 
denote by $C_b(\mathcal{X})$ the linear space of all continuous bounded 
functions $f:\mathcal{X}\to\mathbb{R}$. 
We consider the \emph{superhedging functional}
\begin{equation}\label{superhedgingneu}
\phi(f):=\inf\Big\{\int h\,d\mu_0: h\ge f\mbox{ for some 
}h\in\mathcal{H}\Big\}
\end{equation}
for $f\in C_b(\mathcal{X})$, where $\mu_0\in\mathcal{P}(\mathcal{X})$ is a pricing measure and $\mathcal{H}\subseteq C_b(\mathcal{X})$. Throughout this section we assume that $\mathcal{H}$ is a linear space which contains the constants (i.e.~the constant functions). 

In order to derive a dual representation, we assume that $\phi$ is continuous from above, i.e.~$\phi(f_n)\downarrow 0$ for every sequence $(f_n)$ in $C_b(\mathcal{X})$ such that $f^n\downarrow 0$.
By the nonlinear Daniell-Stone theorem it has a representation 
\begin{equation}\label{superhedgingdualneu}
\phi(f)=\max_{\mu\in\mathcal{Q}} \int f\,d\mu
\end{equation}
for all $f\in C_b(\mathcal{X})$, and the nonempty set $\mathcal{Q}=\big\{\mu\in\mathcal{P}(\mathcal{X}):\int h\,d\mu=\int h\,d\mu_0\mbox{ for all }h\in\mathcal{H}\big\}$. In particular $\mu_0 \in \mathcal{Q}$. The problems \eqref{superhedgingdualneu} and \eqref{superhedgingneu} are in duality and we refer to \eqref{superhedgingdualneu} as the primal and \eqref{superhedgingneu} as the dual formulation. For the details we refer to the Appendix \ref{sec:DS}. There it is outlined how the duality extends to unbounded functions. However, for the sake of readability we focus on $C_b(\mathcal{X})$.

The following example illustrates the basic setting:		 
\begin{example}\label{ex:transport}
	Let $\mathcal{X} = \mathbb{R}^d$, and denote by $\Pi(\mu_1, ..., 
	\mu_d)$ the set of all $\mu\in\mathcal{P}(\mathbb{R}^d)$ with first marginal $\mu_1$, 
	second marginal $\mu_2$, etc. In the following examples, under the 
	assumption that $\mathcal{Q}\neq\emptyset$ it is straightforward to 
	verify that the corresponding superhedging functional is continuous from above. 
	\begin{itemize}
		\item[(a)] (Multi-marginal) optimal transport 
		\cite{kantorovich1942translocation,villani2008optimal}:
		\begin{align*}\mathcal{Q} =\, &\Pi(\mu_1, ..., \mu_d),\\
		\mathcal{H} =\, &\{ h \in C_b(\mathbb{R}^d) : h(x_1, ..., x_d) = h_1(x_1) + ... + 
		h_d(x_d) \text{ for all } (x_1, ..., x_d) \in \mathbb{R}^d\\ 
		&\text{ and some } h_i \in C_b(\mathbb{R})\}\end{align*}
		\item[(b)] Martingale optimal transport 
		\cite{beiglbock2013model,galichon2014stochastic}:		
		\begin{align*}
		\mathcal{Q} =\, &\{ \mu \in \Pi(\mu_1, ..., \mu_d) : \text{the 
			canonical process on $\mathbb{R}^d$ is a $\mu$-martingale}\}\\ 
		\mathcal{H} =\, &\{ h \in C_{\kappa}(\mathbb{R}^d) : h(x_1, ..., x_d) = \sum_{i=1}^d 
		h_i(x_i) + \sum_{i=2}^{d} 
		g_{i}(x_1,\dots,x_{i-1})\cdot(x_{i}-x_{i-1})\\ &\text{ for all } 
		(x_1, ..., x_d) \in \mathbb{R}^d \text{ and } \text{some } h_i \in 
		C_b(\mathbb{R})\text{ and } g_{i}\in 
		C_b(\mathbb{R}^{i-1})\}\end{align*}
		where $C_\kappa(\mathbb{R}^d)$ denotes the space of all continuous functions of linear growth corresponding to $\kappa(x):=1+|x|$, see Appendix \ref{sec:DS}. By Strassen's 
		theorem \cite{strassen1965existence} the set $\mathcal{Q}$ is 
		nonempty if the marginals $\mu_1,\dots,\mu_d$ are in convex order.
		
		\item[(c)] Optimal transport with additional constraints: 
		\begin{align*}
		\mathcal{Q} =\, &\{ \mu \in \Pi(\mu_1, ..., \mu_d) : \int g_j 
		\,d\mu = c_j \text{ for all } j = 1,...,N\} \\
		\mathcal{H} =\, &\{ h \in C_b(\mathbb{R}^d) : h(x_1,..., x_d) = \sum_{i=1}^d 
		h_i(x_i) + \sum_{j=1}^N \lambda_j \left( g_j(x_1, ..., x_d) - 
		c_j\right) \\ &\text{ for some } h_i \in C_b(\mathbb{R}), 
		\lambda_j \in \mathbb{R}\} 
		\end{align*}
		for some $g_1,\dots,g_N \in C_b(\mathbb{R}^d)$ and $c_1,\dots,c_N \in 
		\mathbb{R}$. For related problems we refer to 
		\cite{bartl2017sharpness} and the references therein.		
	\end{itemize}
\end{example}

\subsection{Regularization of the superhedging functional by penalization}
Our goal is to regularize the superhedging functional $\phi$ by considering the convolution
\begin{align}
\phi_{\theta,\gamma}(f)&:=\inf_{h \in C_b(\mathcal{X})}\big\{\phi(h)+\psi_{\theta,\gamma}(f-h)\big\}\nonumber \\
&=\inf_{h\in\mathcal{H}}\Big\{\int h\,d\mu_0+\int \beta_\gamma(f-h)\,d\theta\Big\}\label{defphithetagamma}
\end{align}
where $\psi_{\theta,\gamma}(f):=\int \beta_\gamma(f)\,d\theta$ for a \emph{sampling measure} $\theta\in\mathcal{P}(\mathcal{X})$, and
$\beta_\gamma(x):=\frac{1}{\gamma}\beta(\gamma x)$
is a \emph{penalty function}  which is parametrized by $\gamma>0$. We assume that $\beta\colon\mathbb{R}\to\mathbb{R}_+$ is a
differentiable nondecreasing convex function   such that 
$\lim_{x\to\infty} \beta(x)/x =\infty$. 
Its convex conjugate 
\[\beta^\ast_\gamma(y):=\sup_{x\in\mathbb{R}} 
\{xy-\beta_\gamma(x)\}\quad \mbox{for all }y\in\mathbb{R}_+,\] satisfies $\beta^\ast_\gamma(y)=\beta^\ast(y)/\gamma$.
Common examples are
\begin{itemize}
	\item[(a)] the \textsl{exponential penalty function} $\beta(x)=\exp(x-1)$ with conjugate  $\beta^\ast(y)=y\log(y)$,
	\item[(b)] the \textsl{$L^p$ penalty function} $\beta(x)=\frac{1}{p}(\max\{0,x\})^p$ with conjugate  $\beta^\ast(y)=\frac{1}{q}y^q$ where $q=\frac{p}{p-1}$ for some $p>1$.  
\end{itemize}
In case that $\mathcal{H}=\mathbb{R}$ the functional \eqref{defphithetagamma} is a 
so-called \emph{optimized certainty equivalent}, see Ben-Tal and Teboulle  \cite{ben2007old}.
In the following result we show the 
dual representation of the regularized superhedging functional $\phi_{\theta,\gamma}$ and its convergence to $\phi$.

\begin{theorem}\label{thm:reg}
	Let $f\in C_b(\mathcal{X})$. Suppose there exists $\pi\in\mathcal{Q}$ 
	such that $\pi\ll\theta$ and $\int \beta^\ast\big(\frac{d\pi}{d\theta}\big)\,d\theta<\infty$. Then
	\begin{equation}\label{eq:dualreg} 
	\phi_{\theta,\gamma}(f)=\max_{\mu\in\mathcal{Q}}\Big\{\int f\,d\mu-\frac{1}{\gamma}\int \beta^\ast\Big(\frac{d\mu}{d\theta}\Big)\,d\theta\Big\} .
	\end{equation} 
	Moreover, 
	\begin{equation}\label{eq:boundmain} \phi_{\theta,\gamma}(f)-\frac{\beta(0)}{\gamma}\le \phi(f)\le \phi_{\theta,\gamma}(f)+\frac{1}{\gamma}\int\beta^\ast\Big(\frac{d\mu_\varepsilon}{d\theta}\Big)\,d\theta+\varepsilon
	\end{equation}
	whenever $\mu_\varepsilon\in\mathcal{Q}$ is an $\varepsilon$-optimizer of \eqref{superhedgingdualneu} such that $\mu_\varepsilon\ll\theta$ and $\int \beta_\gamma^\ast \big(\frac{d\mu_\varepsilon}{d\theta}\big) d\theta < \infty$. 
	
	If $\hat h\in\mathcal{H}$ is a minimizer of \eqref{defphithetagamma} then $\hat{\mu} \in \mathcal{P}(\mathcal{X})$ defined by
	\begin{equation}\label{dualoptimizer}
	\frac{d\hat\mu}{d\theta}:=\beta_{\gamma}^\prime(f-\hat h)
	\end{equation}
	is a maximizer of \eqref{eq:dualreg}.
\end{theorem}

\subsection{Approximation of the superhedging functional}
In this subsection we consider a sequence $\mathcal{H}^1\subseteq \mathcal{H}^2\subseteq 
\cdots$ of subsets of $\mathcal{H}$, and 
set $\mathcal{H}^\infty:=\bigcup_{m\in\mathbb{N}} \mathcal{H}^m$. For each 
$m\in\mathbb{N}\cup\{+\infty\}$,
we define the \textsl{approximated superhedging functional} by
\begin{equation}\label{superhedging:approx2}
\phi^m(f):=\inf\Big\{\int h\,d\mu_0\colon h\ge f\;\mbox{for some 
}h\in\mathcal{H}^m\Big\}.
\end{equation}
For the approximation of $\phi(f)$ by $\phi^m(f)$, we need the 
following density condition on $\mathcal{H}^\infty$.
\begin{description}
	\item[Condition (D):] For every $\varepsilon>0$ and $\mu\in 
	\mathcal{P}(\mathcal{X})$ holds
	\begin{itemize}
		\item[(a)] for every $h\in\mathcal{H}$ there exists 
		$h^\prime\in\mathcal{H}^\infty$ such that $\int 
		|h-h^\prime|\,d\mu\le\varepsilon$,
		\item[(b)] there exists $h^{\prime\prime}\in\mathcal{H}^\infty$ 
		such that $1_{K^c}\le h^{\prime\prime}$ and
		$\int h^{\prime\prime}\,d\mu\le\varepsilon$ for some compact subset 
		$K$ of $\mathcal{X}$.
	\end{itemize}
\end{description}
In Section \ref{sec3} we will discuss Condition (D) in the context of 
multilayer feedforward networks. The condition allows for the following 
approximation result.
\begin{proposition}\label{prop:approx}
	Assume that $\mathcal{H}^{\infty}$ is a linear space which contains the constants.
	Under Condition (D) one has 
	\[\lim_{m\to\infty}\phi^m(f)=\phi^\infty(f)=\phi(f)\]
	for all $f\in C_b(\mathcal{X})$.
\end{proposition}
Given a sampling measure $\theta$ and a parametrized penalty function
$\beta_\gamma$ as in the previous subsection,
we define the approximated version of the 
regularized superhedging functional by 
\begin{equation}\label{eq:hedg:app:reg}
\phi^m_{\theta,\gamma}(f)=\inf_{h\in\mathcal{H}^m}\Big\{\int h\,d\mu_0+\int \beta_\gamma(f-h)\,d\theta\Big\}
\end{equation}
for all $f\in C_b(\mathcal{X})$. As a consequence of the two approximative steps $\phi_{\theta, \gamma}(f) \rightarrow \phi(f)$ for $\gamma \rightarrow \infty$ in Theorem \ref{thm:reg} and $\phi^m(f) \rightarrow \phi(f)$ for $m \rightarrow \infty$ in Proposition \ref{prop:approx} we get the following convergence result. 

\begin{proposition}
	\label{prop:unif}
	Suppose that $\mathcal{H}^\infty$ satisfies Condition (D) and for every $\varepsilon>0$ there exists an $\varepsilon$-optimizer $\mu_\varepsilon$ of  \eqref{eq:dualreg} such that $\mu_\varepsilon\ll \theta$ and  $\int \beta^\ast\big(\frac{d\mu_\varepsilon}{d\theta}\big)\,d\theta<+\infty$. Then, for every $f\in C_b(\mathcal{X})$ one has
	$\phi^m_{\theta, \gamma}(f) \rightarrow \phi(f)$ for $\min\{m, \gamma\} \rightarrow \infty$.
\end{proposition}

The existence of such $\varepsilon$-optimizers as required in Theorem \ref{thm:reg} and Proposition \ref{prop:unif} is for example established in \cite{bindini2019smoothing} in the context of multi-marginal optimal transport problems in $\mathbb{R}^d$ with absolutely continuous marginals. In general, the existence of such $\varepsilon$-optimizers crucially depends on the choice of $\theta$, see also Example \ref{example_nouniform} for a simple illustration.

\section{Modelling finite dimensional subspaces with multilayer feedforward 
	networks}\label{sec3}
This section explains the specific choice of approximative subspaces as built by neural 
networks.
Generally, a feasible 
alternative to neural networks is to build these spaces via basis functions, 
like polynomials, 
which is for example pursued in \cite{henry2013automated} in the context of 
martingale optimal transport. In contrast to a 
basis approach, where functions are represented as a weighted sum over fixed 
basis functions, neural networks rely on the composition of layers of simple 
functions. This has shown to be an efficient way to approximate a large class 
of functions with relatively few parameters.
Before going into the results, we give the required 
notation for 
neural networks.

\subsection{Notation} 
The type of neural networks we consider are fully connected feed-forward neural networks. Those are mappings of the form
\[
\mathbb{R}^d \ni x \mapsto A_l \circ \underbrace{\varphi \circ 
		A_{l-1}}_{(l-1).~\text{layer}} 
	\circ ... \circ \underbrace{\varphi \circ A_0}_{\text{1.~layer}} 
	(x)
\]
where $A_i$ are affine transformations and $\varphi : \mathbb{R} \rightarrow \mathbb{R}$ is a nonlinear \textsl{activation} function that is applied elementwise, i.e.~$\varphi((x_1, ..., x_n)) = (\varphi(x_1), ..., \varphi(x_n))$ for $(x_1, ..., x_n) \in \mathbb{R}^n$. 

Regarding dimensions, there is an input dimension $d \in \mathbb{N}$ and a hidden dimension $m \in \mathbb{N}$. This means $A_0$ maps from $\mathbb{R}^d$ to $\mathbb{R}^m$, $A_1, ..., A_{l-1}$ map from $\mathbb{R}^m $ to $\mathbb{R}^m$, and $A_l$ maps from $\mathbb{R}^m$ to $\mathbb{R}$. 
Each affine transformation $A_j$ can trivially be represented as $A_j(x) = M_j x + b_j$ for a matrix $M_j$ and a vector $b_j$. All these matrices and vectors together are the parameters of the network, which can be regarded as an element of $\mathbb{R}^D$ for some $D \in \mathbb{N}$.

We will require the sets which contain all feed-forward neural networks with fixed structure (i.e.~fixed number of layers and fixed dimensions) but unspecified parameter values. We denote by $\Xi \subset \mathbb{R}^D$ the sets of possible parameters for a fixed network structure (where formally, $D$ depends on the structure of the network), and by $N_{l, d, m}(\xi) = A_l \circ \varphi \circ A_{l-1} \circ ... \circ \varphi \circ A_0$ a particular neural network with $l$ layers, input dimension $d$, hidden dimension $m$ and parameters $\xi \in \Xi$. We denote the set of all such networks $N_{l, d, m}(\xi)$ for $\xi \in \Xi$ by $\mathfrak{N}_{l,d, m}(\Xi)$.

In the remainder of this section, we work with a fixed number of layers and input dimension, but allow for growing hidden dimension.
For different hidden dimensions $m$, denote by $\Xi_m$ the corresponding parameter sets. We define
\[
\mathfrak{N}_{l, d} := \bigcup_{m\in\mathbb{N}} \mathfrak{N}_{l,d, m}(\Xi_m).
\]
We want this definition to be independent of the precise choices of the parameter sets, which is why we make the standing assumption that the sets $\mathfrak{N}_{l,d,m}(\Xi_m)$ are growing in $m$. One way to make this explicit is: 
\begin{assumption}
	\label{nn_assumption}
	For any $l, d \in \mathbb{N}$ and a sequence of parameter sets $\Xi_1, \Xi_2, ...$, where $\Xi_m$ is regarded as a subset of  $\mathbb{R}^{D_m}$ for some $D_m \in \mathbb{N}$, we will always assume that
	$[-m, m]^{D_m} \subseteq \Xi_m$ and $\mathfrak{N}_{l,d, m}(\Xi_m) \subset \mathfrak{N}_{l,d,{m+1}}(\Xi_{m+1})$ for all $m \in \mathbb{N}$.
\end{assumption}
The only reason why we do not just set $\Xi_m \equiv \mathbb{R}^{D_m}$ is 
that in Proposition \ref{nouniform} we make the assumption of compact parameter sets.
Further, we assume 
\begin{assumption}
	\label{activation_assumption}
	The activation function $\varphi$ is continuous, nondecreasing and satisfies the limit properties
	$\lim_{x\to-\infty}\varphi(x)=0$ 
	and 
	$\lim_{x\to+\infty}\varphi(x)=1$.
\end{assumption}

\subsection{Modelling $\mathcal{H}^m$ via neural networks}\label{modelHm}
In the following we assume that $\mathcal{H}$ is of the form
\[\mathcal{H}=\Big\{\sum_{j=1}^J e_j  h_j\circ \pi_j + a : h_j\in 
C_b(\mathbb{R}^{d_j}), a\in\mathbb{R} \Big\},\]
where $e_j\in C_b(\mathcal{X})$ and  $\pi_j:\mathcal{X} \to \mathbb{R}^{d_j}$ are continuous 
functions for all $j=1,\dots,J$. This form of $\mathcal{H}$ includes many 
different problems, for instance the ones considered in Example 
\ref{ex:transport} (e.g.~in (a) one has 
$\mathcal{H}=\{\sum_{j=1}^d h_j\circ {\rm pr}_k:h_j\in C_b(\mathbb{R})\}$ where 
${\rm pr}_j(x):=x_j$ denotes the projection on the $j$-th marginal component).

We approximate $\mathcal{H}$ by
\[\mathcal{H}^\infty=\Big\{\sum_{j=1}^J e_j  h_j\circ \pi_j + a : h_j\in 
\mathfrak{N}_{l_j, d_j}, a\in\mathbb{R} \Big\},\] and its 
subspaces 
\[\mathcal{H}^m=\Big\{\sum_{j=1}^J 
e_j  h_j\circ \pi_j + a : h_j\in \mathfrak{N}_{l_j,d_j, m}(\Xi_{j, m}), 
a\in\mathbb{R}\Big\}.\]
In this context the problems $\phi^m_{\theta, \gamma}(f)$ are given by
\begin{align*}
\phi^m_{\theta, \gamma}(f)&=\inf_{h\in \mathcal{H}^m}\Big\{\int h \,d\mu_0 + \int 
\beta_\gamma(f-h)\,d\theta  
 \Big\} \\
&=\inf_{a\in\mathbb{R}}\inf_{h_j\in \mathfrak{N}_{l_j,d_j, m}(\Xi_{j, m})}\Big\{\int \sum_{j=1}^J e_j h_j\circ \pi_j 
\,d\mu_0 + 
a + \int \beta_\gamma\Big(f-\sum_{j=1}^J e_j  
h_j\circ\pi_j - a \Big)\,d\theta \Big\} \\
&=\inf_{a\in \mathbb{R}} \inf_{\xi_j \in \Xi_{j, m}} \Big\{\int \sum_{j=1}^J e_j 
N_{l_j, d_j, m}(\xi_j)\circ \pi_j \,d\mu_0 + 
a + \int \beta_\gamma\Big(f-\sum_{j=1}^J e_j  
N_{l_j, d_j, m}(\xi_j)\circ\pi_j - a \Big)\,d\theta \Big\}
\end{align*}
for all $f\in C_b(\mathcal{X})$. The final formulation illustrates that the problem $\phi^m_{\theta, \gamma}(f)$ is now reduced to a finite dimensional problem of finding the optimal parameters in a neural network. Further, the overall objective depends smoothly on the parameters, and the parameters are unconstrained. In short, problem $\phi^m_{\theta, \gamma}(f)$ fits into the framework of machine learning problems that can be numerically solved by standard stochastic gradient descent based methods.

Under the standing Assumptions \ref{nn_assumption} and \ref{activation_assumption}, the following lemma establishes situations when Condition (D), which is required for Proposition \ref{prop:approx}, is satisfied in the neural network setting.
\begin{lemma}
	\label{lem:CondD}
	\item[(a)] $\mathcal{H}^{\infty}$ satisfies the first part of 
	Condition $(D)$.
	\item[(b)] If $\mathcal{X} = \mathbb{R}^d = 
	\mathbb{R}^{d_1} \times ... \times \mathbb{R}^{d_{J_0}}$ and $\pi_j = {\rm pr}_j,~ 
	e_j 
	= 1$ for $j=1,..., J_0 \leq J$, where ${\rm pr}_j$ is the projection 
	from 
	$\mathbb{R}^d$ to $j$-th marginal component 
	$\mathbb{R}^{d_j}$, then $\mathcal{H}^{\infty}$ satisfies the second part of 
	Condition $(D)$. Further, the second part of Condition $(D)$ is trivially satisfied whenever $\mathcal{X}$ is compact.
\end{lemma}
Notably, part (b) can be seen as a large, but still exemplary case. Intuitively, the second part of Condition (D) is satisfied whenever the space $\mathcal{H}^\infty$ is rich enough.
\begin{remark}
	Later in the numerics we will usually work with a ReLU activation function, 
	i.e.~$\varphi(x) = \max\{0, x\}$. While this does not satisfy the 
	latter limit property of Assumption \ref{activation_assumption}, this is easily 
	amendable: 
	Basically, throughout the whole theory the assumptions will only be used to guarantee existence of neural networks with certain properties.
	Given Assumption \ref{activation_assumption}, we will only require two layers ($l=1$) to obtain the necessary results.
	In the numerics however, we use more layers. If more layers are given, one can also bundle several layers and regard them as one layer, with a different activation function. For example:
	\[
	A_l \circ \underbrace{\varphi \circ 
		A_{l-1}
	\circ ... \circ A_1 \circ \varphi}_{\overline{\varphi}} \circ A_0 	
	\]
	Whenever $\overline{\varphi}$ is a mapping of the form $(x_1, ..., x_m) \mapsto (\overline{\varphi}(x_1), ..., \overline{\varphi}(x_m))$, an $(l+1)$-layer network with activation function $\varphi$ can represent any function that a two layer network with activation function $\overline{\varphi}$ can represent. For $\varphi(x) = \max\{0, x\}$ one can easily see that $\overline{\varphi}(x) = \min\{1, \max\{0, x\}\}$ is feasible, which satisfies Assumption \ref{activation_assumption}.
\end{remark}

\subsection{Convergence}
In this section we study in what sense $\phi^m_{\theta, \gamma}(f)$ converges to $\phi(f)$ for the approximation by neural networks.

First, we study the case of uniform convergence in $m$ and $\gamma$, i.e.~conditions for the convergence $\phi^m_{\theta, \gamma}(f) \rightarrow \phi(f)$ for $\min\{m, \gamma\} \rightarrow \infty$. This is subject of Remark \ref{uniform} below, which is a summary of results established in Section \ref{sec2} and Section \ref{modelHm}. The two approximative steps leading to uniform convergence are $\phi_{\theta, \gamma}(f) \rightarrow \phi(f)$ for $\gamma \rightarrow \infty$ and $\phi^m(f) \rightarrow \phi(f)$ for $m \rightarrow \infty$.

On the other hand, sometimes the convergence 
$\phi_{\theta, \gamma}(f) \rightarrow \phi(f)$ for $\gamma \rightarrow \infty$ is not satisfied even though practically one obtains a good approximation. One 
such case is given in Example \ref{example_nouniform}. Even if uniform 
convergence does not hold, one can still often connect problems 
$\phi^m_{\theta,\gamma}(f)$ and $\phi(f)$. This is done by the approximative 
steps $\phi^m(f) \rightarrow \phi(f)$ for $m\rightarrow \infty$ and 
$\phi_{\theta,\gamma}^m(f) \rightarrow \phi^m(f)$ for $\gamma \rightarrow 
\infty$, where the latter is subject of Proposition \ref{nouniform}. Here, instead of 
the strong assumption required for $\phi_{\theta, \gamma}(f) \rightarrow 
\phi(f)$, the convergence $\phi_{\theta,\gamma}^m(f) \rightarrow \phi^m(f)$ can 
be shown by assuming that all parameter sets of the neural networks are compact.

\begin{remark}
	\label{uniform}
	Under the assumptions of Lemma \ref{lem:CondD}, Proposition \ref{prop:approx} implies $\phi^m(f) \rightarrow \phi(f)$ for $m\rightarrow \infty$.
	Further, given the existence of $\varepsilon$-optimizers for every $\varepsilon > 0$ as required in Theorem \ref{thm:reg}, convergence $\phi_{\theta, \gamma}(f) \rightarrow \phi(f)$ for $\gamma \rightarrow \infty$ holds.
	Given both assumptions, Proposition \ref{prop:unif} yields $\phi_{\theta, \gamma}^m(f) 
	\rightarrow \phi(f)$ for $\min\{m, \gamma\} \rightarrow \infty$. The convergence $\phi_{\theta, \gamma}^m(f) \rightarrow \phi_{\theta, \gamma}(f)$ for $m\rightarrow \infty$ is a trivial consequence.
\end{remark}

\begin{example}
	\label{example_nouniform}
	Let $\mathcal{X} = [0, 1]^2$, $\mu_1 = \mu_2 = \delta_0$ and $f(x_1, x_2) = -|x_1 - x_2|$. Let $\mathcal{Q} = \Pi(\mu_1, \mu_2)$ be the set of all measures in $\mathcal{X}$ with first marginal $\mu_1$ and second marginal $\mu_2$, so that 
	\[
	\phi(f) = \sup_{\mu \in \Pi(\mu_1, \mu_2)} \int f d\mu
	\]
	Obviously, $\phi(f) = f(0, 0) = 0$. Note that $\mathcal{Q} = \{\mu_1 \otimes \mu_2 \}$ so that $\mu_0 = \mu_1 \otimes \mu_2 = \delta_{(0, 0)}$.
	
	Consider two possible reference measures, $\theta^{(1)} = \mathcal{U}([0, 1]^2)$ being the uniform distribution on $[0, 1]^2$, and $\theta^{(2)} = \mu_1 \otimes \mu_2 = \delta_{(0, 0)}$.
	For $\theta^{(2)}$ it is obvious that the existence of $\varepsilon$-optimizers as required in Theorem \ref{thm:reg} is given, since $\theta^{(2)}$ itself is the optimizer of $\phi(f)$. Hence $\phi_{\theta^{(2)}, \gamma}(f) \rightarrow \phi(f)$ for $\gamma \rightarrow \infty$ holds.
	
	On the other hand, there does not exist $\nu \in \Pi(\mu_1, \mu_2)$ with $\nu \ll \theta^{(1)}$, and hence $\phi_{\theta^{(1)}, \gamma}(f) = -\infty$.
	However, by first approximating $\phi(f)$ by $\phi^m(f)$, the functional becomes smoother: Roughly speaking,
	the marginal constraints are slightly relaxed. This becomes obvious when 
	studying the dual formulations
	\begin{align*}
	\phi_{\theta^{(1)},\gamma}(f) &=  \inf_{h_1, h_2 \in C_b([0, 1])} \Big\{ h_1(0) + h_2(0) + \int \beta_\gamma(f-h_1-h_2) \,d\theta^{(1)} \Big\}\\
	\phi^m_{\theta^{(1)},\gamma}(f) &= \inf_{h_1, h_2 \in \mathfrak{N}_{l, 1, m}(\Xi_m)} \Big\{ h_1(0) + h_2(0) + \int \beta_\gamma(f-h_1-h_2) \,d\theta^{(1)}\Big\}
	\end{align*}
	While one easily finds sequences  of functions in $C_b([0, 1])$ so that the 
	values at $0$ go to minus infinity but the penalty term stays bounded, this 
	is impossible with functions in $\mathfrak{N}_{l, 1, m}(\Xi_m)$ given 
	that the activation function is continuous and the parameter sets are 
	compact. So there is hope to establish the convergence 
	$\phi^m_{\theta^{(1)},\gamma}(f) \rightarrow \phi^m(f)$ for $\gamma \rightarrow 
	\infty$, which will indeed be a consequence of the following result.
\end{example}

\begin{proposition}
	\label{nouniform}
	Fix $m \in \mathbb{N}$. Given that all parameter sets $\Xi_{j, m}$ for $j=1, ..., J$ of the neural networks occurring in $\mathcal{H}^m$ are compact and $\theta$ is strictly positive (i.e.~$\theta$ gives positive mass to every non-empty open set), it holds $\phi^m_{\theta,\gamma}(f) \rightarrow \phi^m(f)$ for $\gamma \rightarrow \infty$.
\end{proposition}

\section{Numerical Examples}
\label{sec:NumEx}
This section aims at showcasing how various frequently studied problems that 
fall into the theoretical framework of the previous sections can be implemented 
simply and effectively with neural networks. The examples focus on toy problems 
that 
allow an objective evaluation of the numerical results and give the reader an 
idea about the strengths and weaknesses of the presented approach. We chose a 
very basic implementation using Tensorflow and the Adam optimizer.\footnote{All 
	used code is available on 
	\url{https://github.com/stephaneckstein/transport-and-related}.} As for the 
	network architecture: In all the 
examples, $\mathcal{H}$ is as described in Section \ref{sec3}, and 
$\mathfrak{N}_{l_k,m}$ always approximates $C_b(\mathbb{R}^d)$. To approximate 
$C_b(\mathbb{R}^d)$ we use a five layer ($l=4$ in the previous chapter) ReLU-network with hidden dimension 
$64 \cdot d$. 
We did not perform a hyper parameter search to obtain this architecture, but 
rather oriented ourselves at papers with comparable settings (e.g.~at 
\cite{weinan2017deep, gulrajani2017improved, seguy2017large}). Notably, increasing the complexity (number of layers or hidden dimension) further did not change the numerical results significantly in the cases tested, so we believe the structure chosen to be adequate for the problems considered.

Simply put, the implementation works as follows: We perform a normal stochastic 
gradient type optimization (outsourced to the Adam optimizer) for a certain 
number of iterations to find near
optimal parameters of the network. At each iteration during this process, the 
expectations in the objective function are replaced by averages over a fixed 
number (called \textsl{batch size}) of random points from the respective 
distributions. To obtain the numerical approximation 
$\hat{\phi}_{\theta,\gamma}^{m}(f)$ of $\phi_{\theta,\gamma}^{m}(f)$, we finally 
average 
the sample objective values over the last roughly $5\%$ of iterations. This is referred to as the dual value.
Alternatively, one can use formula \eqref{dualoptimizer} to obtain sample points from an approximate optimizer $\nu^\ast$ of the primal problem and numerically evaluate $\int f \,d\nu^\ast$, which is referred to as the primal value (more details on how to work with such an approximative optimizer $\nu^\ast$ is given in Section \ref{sec::Rearrangement}).  If not stated otherwise, all reported values are dual values.

The numerical procedure we use can likely be improved by fine-tuning parameters 
or by using 
more complex network architectures. For example batch normalization is applied 
in a related setting in \cite{buhler2018deep} which appears to significantly 
speed up the optimization.

\subsection{Optimal transport and Fr\'echet-Hoeffding bounds}
\label{subsec::Frechet}
With this first problem, we study the effects of different penalty functions, penalty factor, batch size and number of iterations of the Adam optimizer.
\begin{figure}[t]\hspace{-0.85cm}
	\mbox{\includegraphics[width=.55\textwidth]{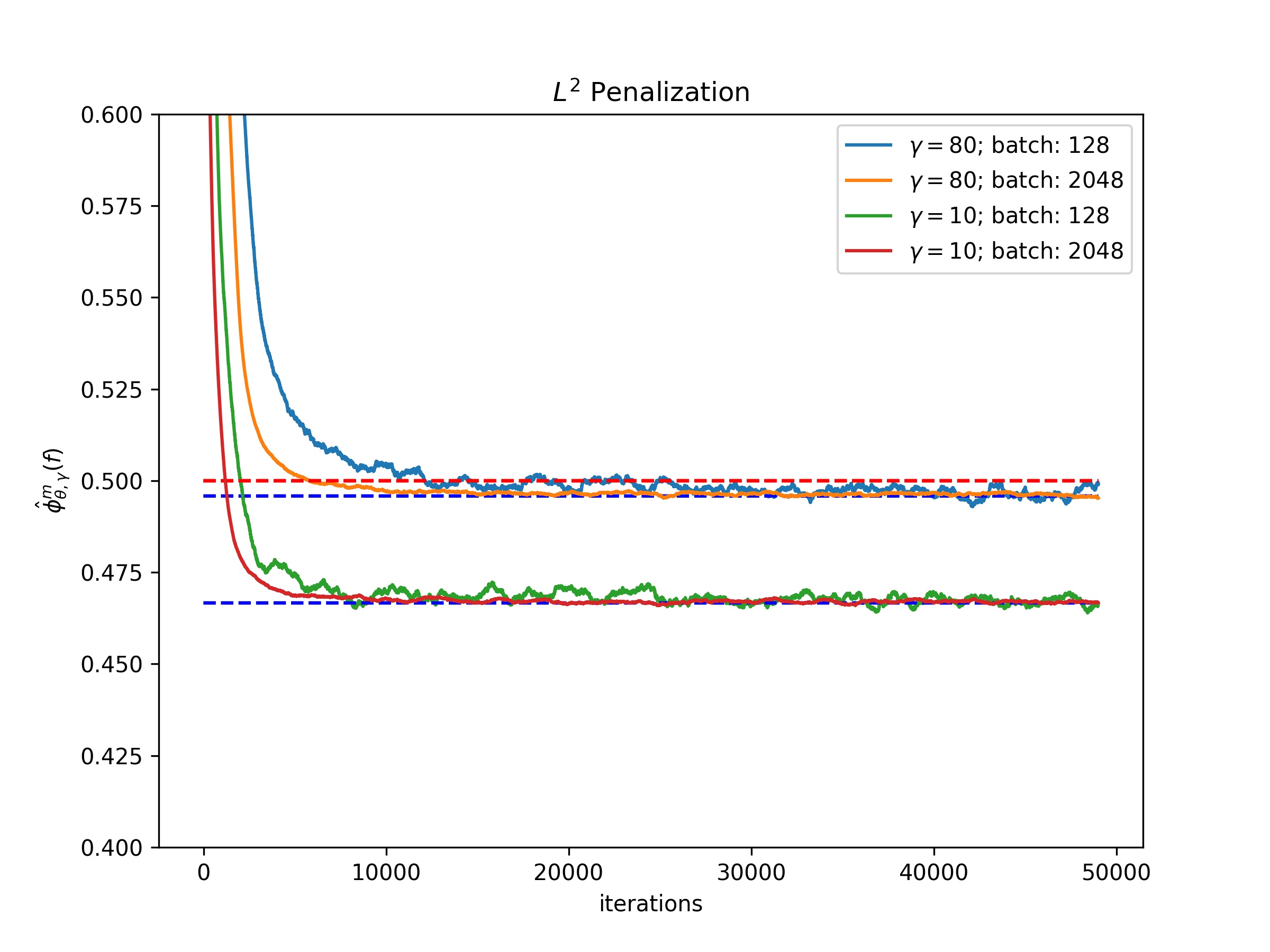}
		\includegraphics[width=.55\textwidth]{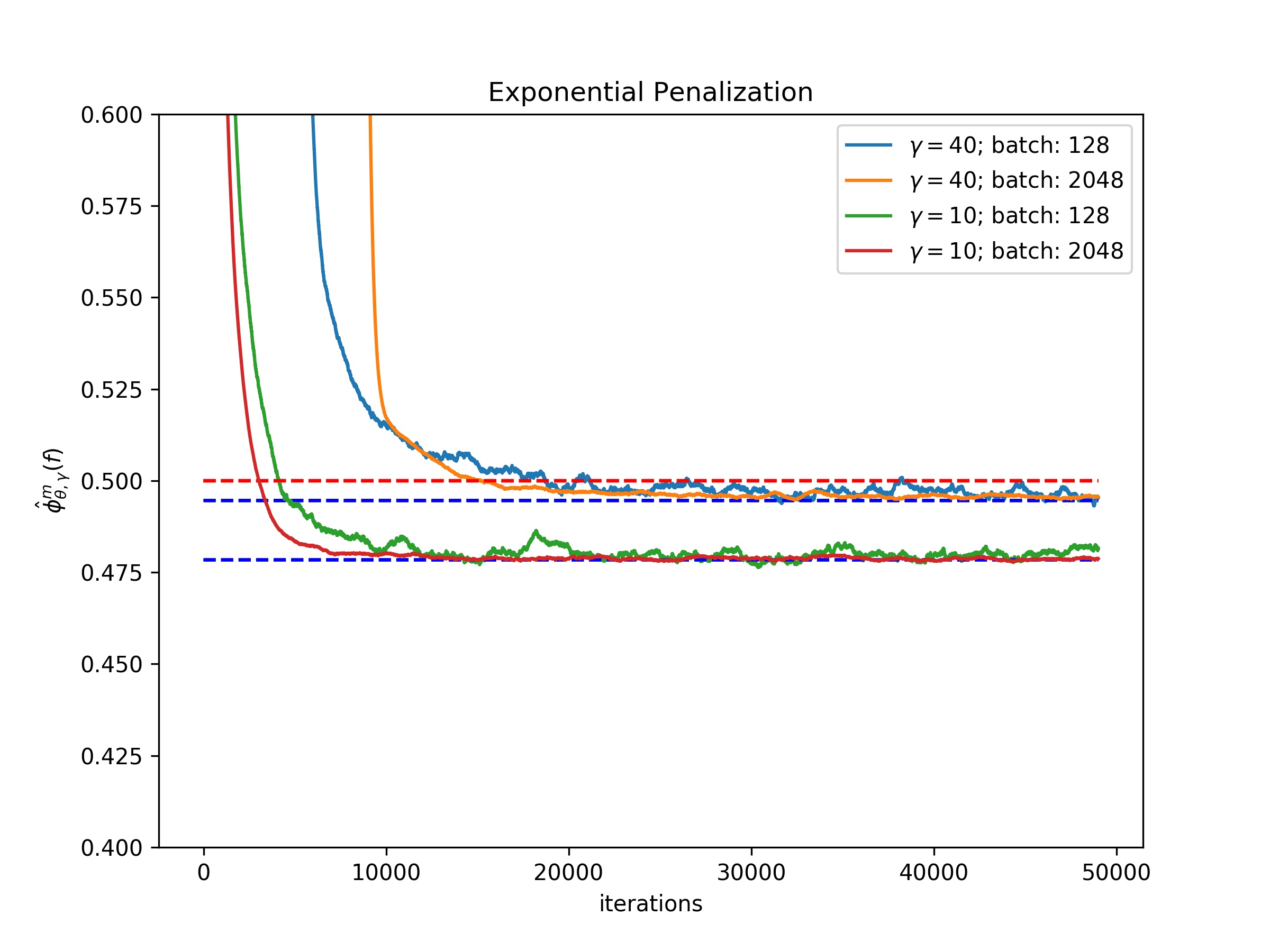}}
	\caption{Fr\'echet-Hoeffding bounds: $d = 2,~z_1 = 0.5,~ z_2 = 0.75$. 
		Comparison of $L^2$ penalty 
		function 
		$\beta_{\gamma}(x) = 
		\gamma \max\{0,x\}^2$ and exponential penalty 
		function $\beta_{\gamma}(x) = 
		\frac{\exp(\gamma x - 1)}{\gamma}$. The values plotted are running 
		averages 
		over 
		the last 1000 iterations. The dotted red line is the true 
		value $\phi(f)$. The dotted blue lines are bounds 
		from 
		below for $\phi_{\theta,\gamma}(f)$ obtained by Equation \eqref{eq:boundmain} in Theorem \ref{thm:reg} for the 
		respective choices of $\gamma$.}
	\label{FigureFH1}
\end{figure}
Let 
$\mathcal{X} = [0,1]^d,~\theta = 
\mathcal{U}\left([0,1]^d\right)$ (where $\mathcal{U}(\cdot)$ 
denotes the uniform distribution) and $\mathcal{Q} = \{\nu \in \mathcal{P}(\mathcal{X}): 
\nu_i 
= 
\mathcal{U}\left([0,1]\right)\}$, where $\nu_i$ is the $i$-th marginal of $\nu$. For some fixed
$z \in [0,1]^d$, define the function $f : [0,1]^d \rightarrow 
\mathbb{R}_+$ by\footnote{The function $f$ here is not continuous. Since the optimal transport problem is continuous from 
	below (see e.g.~\cite{kellerer1984duality}), the representation 
	\eqref{superhedgingdualneu} nevertheless holds for all bounded measurable functions 
	$f$.}
\[
f(x) =  \left\{ \begin{array}{ll} 1,& \text{if } x_i \leq z_i \text{ for all 
}i \in \{1,2,...,d\}, \\ 
0,& \text{else.} \end{array} \right.
\]
The value $\phi(f) = \sup_{\nu \in \mathcal{Q}} \int f \,d\nu$ 
corresponds to the maximum value of a 
$d$-dimensional copula at point $z$. By the Fr\'echet-Hoeffding bounds we have 
an analytical solution to this problem, which is
\[
\phi(f) = \min_{i\in \{1,...,d\}} z_i.
\]

In Figure \ref{FigureFH1} we 
observe how $\hat{\phi}_{\theta,\gamma}^m(f)$ depends on the number of 
iterations of 
the Adam optimizer and the batch size. We observe that 
while higher batch sizes lead to more stable convergence, the speed of 
convergence appears not strongly related to batch size. This suggests that 
increasing batch sizes might lead to both quick and finally stable 
performance.\footnote{See also \cite{smith2017don} and references therein for 
related concepts on how to optimally tune the optimization procedure. In this 
paper however, we decided to stick with standard parameters of the Adam 
optimizer and fixed batch size. This is done to avoid another layer of 
complexity when evaluating the numerical results.} Since $L^2$ penalization 
appears more 
stable, we will mostly use this penalization for the rest of the applications. Further, the figure illustrates that the numerical solutions appear to approximately obtain the lower bounds for $\phi_{\theta, \gamma}(f)$ as given by Equation \eqref{eq:boundmain} in Theorem \ref{thm:reg}. I.e.~one approximately has $\phi(f) \approx \phi_{\theta, \gamma}^m(f) + \frac{1}{\gamma} \int \beta^\ast(\frac{d\hat{\mu}}{d\theta}) d\theta$ where $\hat{\mu}$ is an optimizer of $\phi(f)$.\footnote{Here, $\hat{\mu}$ is chosen as the optimizer which is uniform on the cubes $[0, z]$ and $[z, 1]$.}

\subsection{Multi-marginal optimal transport}
\begin{table}
	\label{tableComparison}
	\centering
	\begin{tabular}{@{}lcccccccccr}
		& & \multicolumn{2}{c}{LP} & \phantom{abc} & \multicolumn{2}{c}{NN} & \phantom{abc} & \multicolumn{1}{c}{RKHS} & \phantom{abc} & \multicolumn{1}{c}{Ref} \\
		\cmidrule{3-4} \cmidrule{6-7} \cmidrule{9-9} \cmidrule{11-11}
		&& MC & quantization & & dual & primal & & Laplace & & Com. \\
		\textbf{$(M, D, K)$} &&&&&&&&&&\\\hline
		&&&&&\multicolumn{2}{c}{\textbf{p = q = 2}}&&&& \\
		\midrule
		$(2, 1, 1)$ & & \begin{tabular}{@{}c@{}}0.403\\ (0.084) \end{tabular} & \begin{tabular}{@{}c@{}}0.408\\ (0.026) \end{tabular} & & 0.413 & 0.401 & & \begin{tabular}{@{}c@{}}0.364\\ (0.006) \end{tabular} & & 0.405 \\\hline
		$(2, 1, 6)$ & & \begin{tabular}{@{}c@{}}3.337\\ (0.320) \end{tabular} & 
		\begin{tabular}{@{}c@{}}3.263\\ (0.115) \end{tabular} & & 3.279 & 3.258 
		& & \begin{tabular}{@{}c@{}}2.444\\ (0.018) \end{tabular} & & 3.269 
		\\\hline
		$(5, 2, 6)$ & & \begin{tabular}{@{}c@{}}8.978\\ (8.233) \end{tabular} & \begin{tabular}{@{}c@{}}3.073\\ (0.231) \end{tabular} & & 3.123 & 3.041 & & \begin{tabular}{@{}c@{}} DNC \end{tabular} & & - \\\hline
		&&&&&\multicolumn{2}{c}{\textbf{p = 1, q = 2}}&&&& \\
		\midrule
		$(2, 1, 6)$ & & \begin{tabular}{@{}c@{}}1.536\\ (0.071) \end{tabular} & \begin{tabular}{@{}c@{}}1.537\\ (0.025) \end{tabular} & & 1.537 & 1.531 & & \begin{tabular}{@{}c@{}}1.471\\ (0.009) \end{tabular} & & 1.533 \\\hline
		$(5, 2, 6)$ & & \begin{tabular}{@{}c@{}}2.845\\ (1.314) \end{tabular} & \begin{tabular}{@{}c@{}}1.741\\ (0.064) \end{tabular} & & 1.753 & 1.740 & & \begin{tabular}{@{}c@{}} DNC \end{tabular} & & - \\\hline
		$(10, 3, 6)$ & & \begin{tabular}{@{}c@{}}10.235\\ (3.576) \end{tabular} & \begin{tabular}{@{}c@{}}6.744\\ (0.074) \end{tabular} & & 6.759 & 6.743 & & \begin{tabular}{@{}c@{}} DNC \end{tabular} & & - \\\hline&&&&&&&&&&\\[-1em]
		&&&\multicolumn{6}{c}{\boldmath{$\tilde{f}(x) = f(x) \cdot \sin(\sum_{i=1}^M x_{i, 1})$}}&& \\[0.1em] \hline
		$(5, 2, 6)$ & & \begin{tabular}{@{}c@{}}16.814\\ (0.893) \end{tabular} & \begin{tabular}{@{}c@{}}17.380\\ (0.043) \end{tabular} & & 18.001 & 17.539 & & \begin{tabular}{@{}c@{}} DNC \end{tabular} & & - \\\hline
		$(10, 3, 6)$ & & \begin{tabular}{@{}c@{}}24.618\\ (2.332) \end{tabular} & \begin{tabular}{@{}c@{}}23.615\\ (0.107) \end{tabular} & & 34.235 & 32.521 & & \begin{tabular}{@{}c@{}} DNC \end{tabular} & & - \\\hline&&&&&&&&&&\\[-1em]
	\end{tabular}
	\caption{Multi-marginal optimal transport: Numerical values for $-\phi(f)$ 
		arising from different numerical schemes. The numbers in brackets are 
		empirical standard deviations over 100 runs. LP denotes linear programming, 
		based on either sampling the marginals randomly (MC) or using the 
		quantization approach from \cite[Algorithm 4.5]{pflug2014multistage} to 
		approximate marginals (quantization). The neural network (NN) 
		implementation is based on $L_2$ penalization with $\gamma = M \cdot D 
		\cdot 500$. For the reproducing kernel Hilbert space solution (RKHS) as 
		described in \cite[Algorithm 3]{genevay2016stochastic} we use a Laplace 
		kernel and the same penalization as for 
		the NN-method. For the final two rows, we report numerical values for 
		$-\phi(\tilde{f})$. DNC entries did not converge. The final column are analytical reference values given by the comonotone coupling for two marginals.}
\end{table}
The aim of this example is to compare the approach of this paper with existing methods for a numerically challenging problem.
Let $\mathcal{X} = (\mathbb{R}^D)^M$, where $M$ denotes the number of marginals and $D$ denotes the dimension of each marginal. Let $\mu_i$ for $i=1, ..., M$ be $K$-mixtures of normal distributions with randomly chosen parameters, and define $\mathcal{Q} = \Pi(\mu_1, ..., \mu_M)$. For $p, q \geq 1$ let
\[
f(x) := -\Big(\sum_{j=1}^D \big|\sum_{i=1}^M (-1)^i x_{i, j} \big|^q\Big)^{p/q},
\]
where we write $x = (x_{i, j}) \in \mathcal{X}$ with $i=1, ..., M$, $j=1, ..., D$. Note that for two marginals, one has $-\phi(f) = W_{p, q}^p(\mu_1, \mu_2)$, where $W_{p, q}$ is the Wasserstein-$p$ distance with $L_q$ norm on $\mathbb{R}^d$.

In Table \ref{tableComparison} we compare optimal values to this problem 
arising from different algorithmic approaches. We compare linear programming 
methods based on discretization of the marginals, the neural network method 
presented in this paper, and a reproducing kernel Hilbert space (RKHS) approach 
as described in \cite[Algorithm 3]{genevay2016stochastic}. For the linear 
programming 
methods, we use a maximal 
number of variables of $10^6$, which was around the boundary so that Gurobi 
\cite{gurobi} was still able to solve the resulting linear program on our computer. 
Regarding 
the RKHS algorithm we have to mention that it is the only method that is not 
building on an established package like Tensorflow or Gurobi. Hence efficiency 
with respect to running time and tuning of hyperparameters might be far from 
optimal for this approach. Notably, switching from exponential penalization as 
used in \cite{genevay2016stochastic} to $L^2$-penalization was already a slight 
improvement. For the precise specifications of each algorithm and the problem setting, we refer to the code on \url{https://github.com/stephaneckstein/OT_Comparison}.

Evaluating the results, we find that the neural network based method and the linear programming method with quantization appear to work best. Surprisingly, even for the case with 10 marginals (where the linear program can only use 4 points to approximate each marginal!), the quantization method achieved a similar value as the neural network method for $-\phi(f)$. We believe the reason is that the function $f$ is very smooth which a quantization approach can exploit. Hence we slightly changed $f$ to $\tilde{f}$ in the final two test cases, which makes the function less regular. These are the only cases where the neural network solution and the quantization method strongly differ. In the final case, the quantization approach still has to approximate each marginal distribution with just 4 points, while the neural network method can use millions of points. From this standpoint, one can place higher trust in the neural network solution, even though we have no analytical reference value to compare it against.

Initially, we included a fourth method based on the approach in this paper, but with a polynomial basis instead of neural networks. This performed very badly however (at least when using a standard monomial basis), and hence we omitted the results in this table.

\subsection{Martingale optimal transport}
In martingale optimal transport, the optimal transport problem is extended by 
imposing a martingale constraint on top of marginal constraints.
Dimensions are regarded as discrete time-steps and the measures in 
$\mathcal{Q}$ are distributions of discrete stochastic processes $(X_t)_{t = 
	1,...,d}$ with fixed marginal distributions as well as the condition that the 
process is a martingale.

Here, we consider a simple example with $d = 2$, where an analytical solution 
is known. This example is taken from \cite{alfonsi2017sampling}. Let $\mathcal{X} := [-1,1] 
\times [-2,2]$, $\theta := \mathcal{U}(\mathcal{X})$ and set
\[
\mathcal{Q} := \left\{\nu = \nu_1\otimes K :  \nu_1 = 
\mathcal{U}([-1,1]),~
\nu_2 = 
\mathcal{U}([-2,2]),~ x = \int_{-2}^{2} y K(x,dy) \text{ holds } 
\nu_1\text{-a.s.}\right\}.
\]
For $f = -|x-y|^{\rho}$ one gets $\phi(f) = -1$ for all $\rho > 2$. We implement this problem with $\rho = 2.3$, where we use
the $L^2$ penalty function for different values of $\gamma$. The results are 
shown in Figure \ref{FigureMOT1}. One can see that while for values of $\gamma$ 
up to around 1280, the behavior of the optimal value is approximately as 
predicted by Equation \eqref{eq:boundmain} in Theorem \ref{thm:reg}, in that the error decreases by roughly a 
factor of two if $\gamma$ is increased by a factor of two. For larger values of 
$\gamma$ however, 
numerical instabilities occur and the optimizer cannot find the true optimum. 
This is indicated by the fact that the value 
$\hat{\phi}_{\theta,\gamma}^{m}(f)$ 
is 
above $\phi(f) = -1$.

\begin{figure}[t]
	\begin{center}
		\mbox{\includegraphics[width=.6\textwidth]{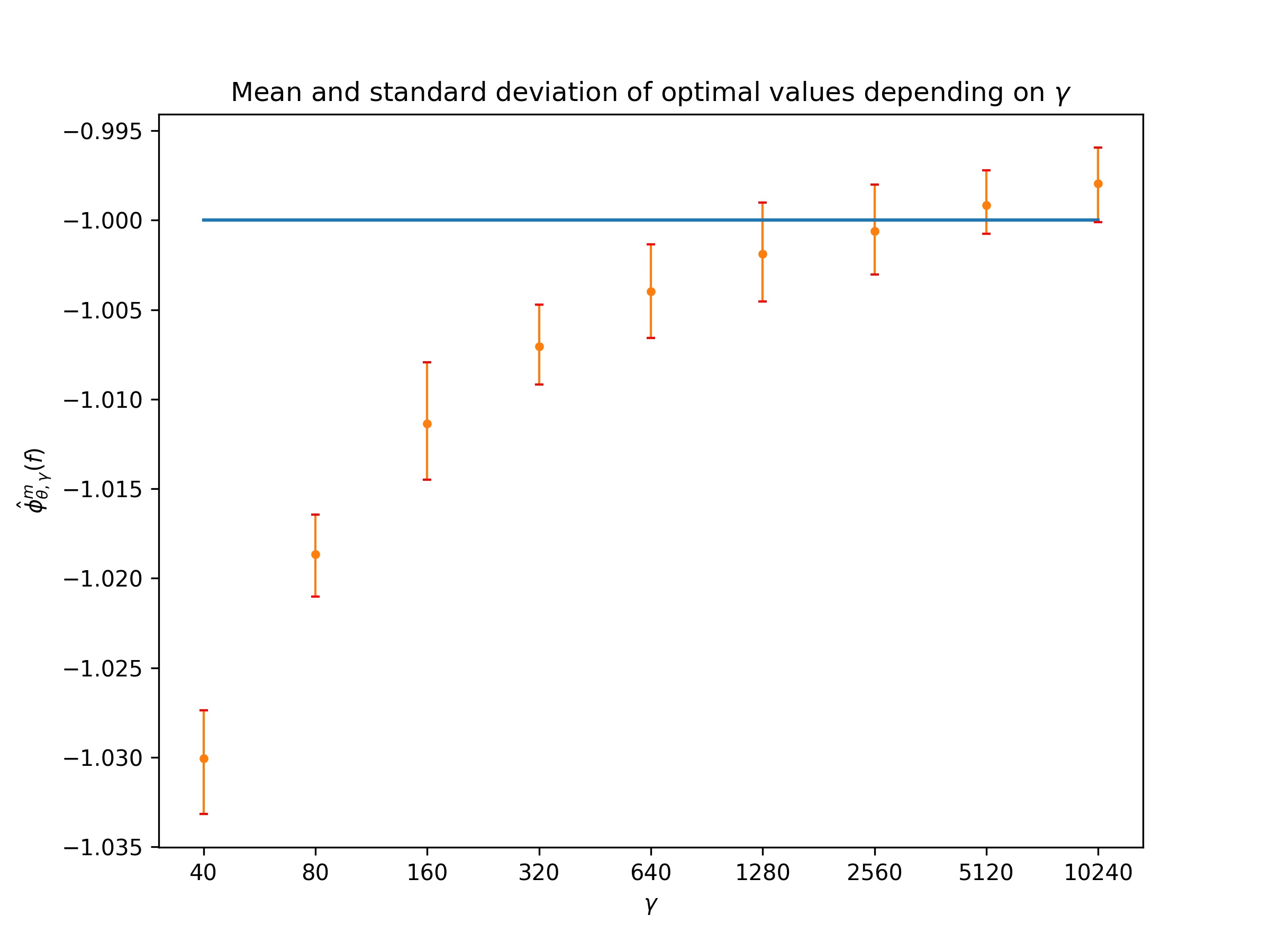}}
		\caption{Martingale optimal transport: Mean numerical optimal values 
			and 
			$95\%$ confidence 
			bounds over 100 independent runs for different values of $\gamma$ 
			($L^2$ penalization). The 
			network is trained for 20000 iterations with batch size 1024. The true 
			optimal value 
			of the unpenalized problem
			is -1.}
		\label{FigureMOT1}
	\end{center}
\end{figure}

\subsection{Portfolio optimization}
\label{subsec::PO}
Consider a market with two assets, where the distribution of returns for each
individual asset is given, but not the joint distribution. An investor wants 
to maximize his or her worst-case utility from investing into the two assets.
Here, the utility of the investor is characterized by a mean-variance 
objective. While the mean is fully characterized by the marginal distributions, the worst case considers all possible variances of the portfolio, which depend on the joint distribution of the assets.

The following example is taken from \cite{pflug2017review}: Let $\mathcal{X} = 
[0,1]\times[0,2]$. Let $\theta_1 = \mathcal{U}([0,1])$ and $\theta_2 = 
\mathcal{U}([0,1]) \circ 
\varphi^{-1}$, where $\varphi(x) = 2 x^2$. Let $\mathcal{Q} = \{\nu \in 
\mathcal{P}(\mathcal{X}): 
\nu_1 = \theta_1,~\nu_2 = \theta_2\}$. We will solve the following robust 
mean-variance portfolio optimization problem
\begin{align*}
\sup_{x \in [0,1]}& - \phi(-f_x) :=\\
\sup_{x \in [0,1]}& \inf_{\nu \in \mathcal{Q}} \int (1-x) \xi_1 + x \xi_2\\ &-
\lambda \left((1-x) \xi_1 + x \xi_2 - (1-x)\int_0^1  \zeta_1 \,\theta_1(d\zeta_1) - 
x\int_0^2
\zeta_2 \theta_2(d\zeta_2)\right)^2 \,\nu(d\xi).
\end{align*}
where $\lambda \geq 0$ is the risk aversion. The integral over the term inside 
the large brackets is the variance of the portfolio. For the analytical 
solution, see Example 1 of \cite{pflug2017review}. We implemented the above 
problem in two ways. For the 
first, we choose the reference measure $\theta^{(1)} = \theta_1 \otimes 
\theta_2$. 
For 
the second, we use the reference measure $\theta^{(2)} = 0.5 \theta^{(1)} + 0.5 
\left(\mathcal{U}([0,1]) \circ (\id, 
\varphi)^{-1}\right)$,\footnote{$\id$ denotes the identity mapping.} i.e.~half 
the product measure, half the perfectly 
correlated measure. The 
second version may correspond to our intuition that the optimal coupling should 
include positive correlation. More precisely: The choice of reference measure 
$\theta$ always has the implicit objective to lead to narrow bounds in Equation \eqref{eq:boundmain} in Theorem \ref{thm:reg}. In this example, if one presumes that an optimal 
measure $\nu^* \in \mathcal{Q}$ has mass near the perfectly correlated 
diagonal, it makes sense to choose a reference measure which puts mass in this 
region, as does $\theta^{(2)}$. The results are reported in Figure 
\ref{FigurePF1}. As expected, the second version yields results closer to the 
analytical solution.
\begin{figure}[t]
	\begin{center}
		\mbox{\includegraphics[width=.5\textwidth]{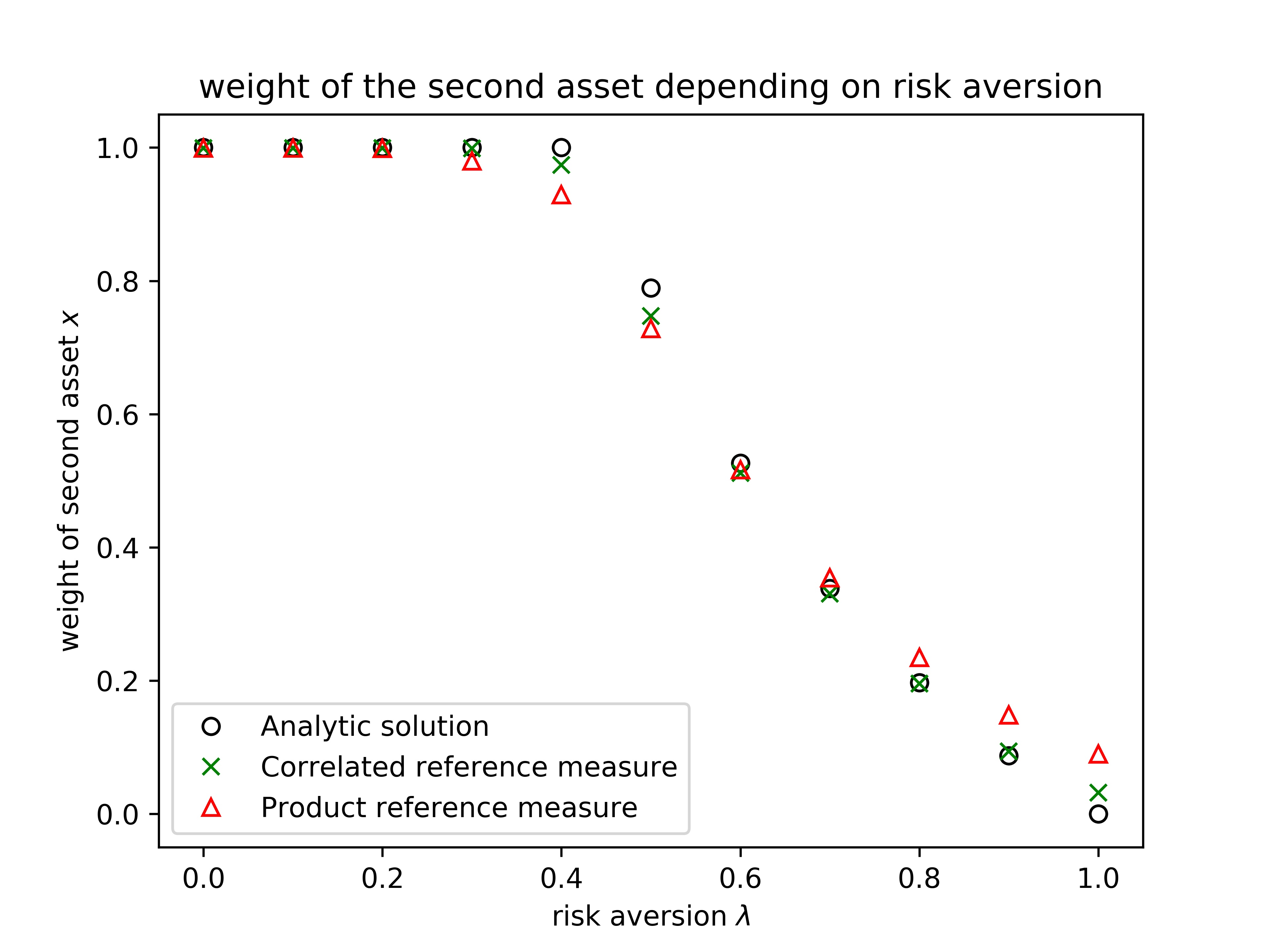}
			\includegraphics[width=.5\textwidth]{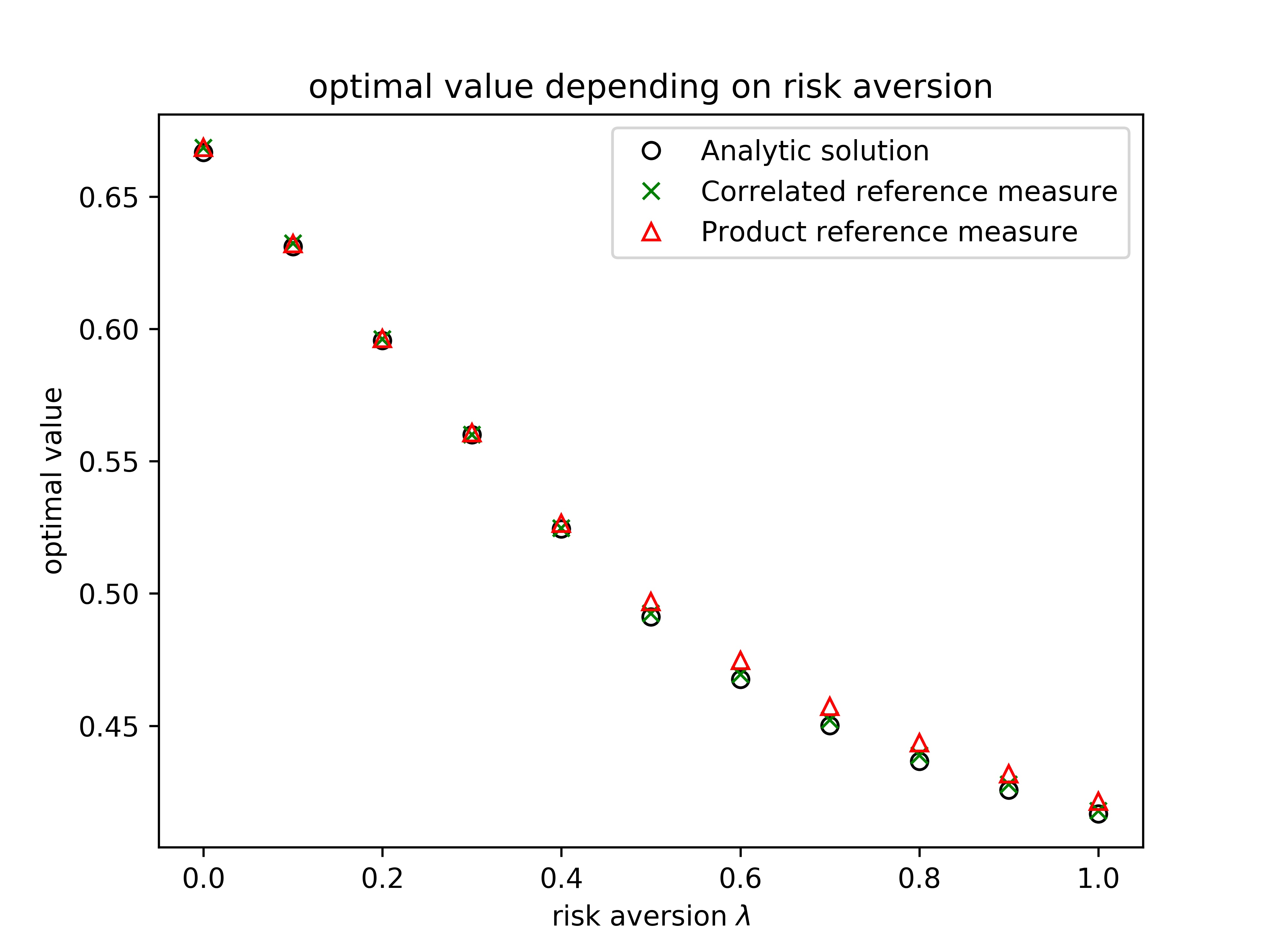}}
		\caption{Portfolio optimization under dependence uncertainty: As 
			reference measure we take either the product measure or a
			positively correlated measure. We use $L^2$ penalization 
			with 
			$\gamma = 160$. The network is trained with batch size $2^{13}$ for 
			40000 iterations.}
		\label{FigurePF1}
	\end{center}
\end{figure}

\subsection{Bounds on the distribution of a sum of dependent random variables}
\label{sec::Rearrangement}
In this section, the objective is to find bounds for the probability 
$\mathbb{P}(X_1 + X_2 + ... + X_d \geq s)$ for some $s\in \mathbb{R}$, where 
the individual distributions of $X_i$ are known, but not their joint 
distribution. This problem is in strong relation to calculating worst- and 
best-case value at risks under dependence uncertainty, see also 
\cite{eckstein2018robust, embrechts2013model, puccetti2012computation}. 

Let $\mathcal{X} = \mathbb{R}^d$. For given marginals $\mu_1, ..., \mu_d \in 
\mathcal{P}(\mathbb{R})$ and $\mathcal{Q} = 
\Pi(\mu_1, ..., \mu_d)$, the problem statement is
\[
\phi(f) := \sup_{\nu \in \mathcal{Q}} \int 1_{\{x_1 + ... + x_d \geq s\}} 
\,\nu(dx_1, ..., dx_d).
\]

\begin{figure}[t]
	\begin{center}
		\mbox{\includegraphics[width=.7\textwidth]{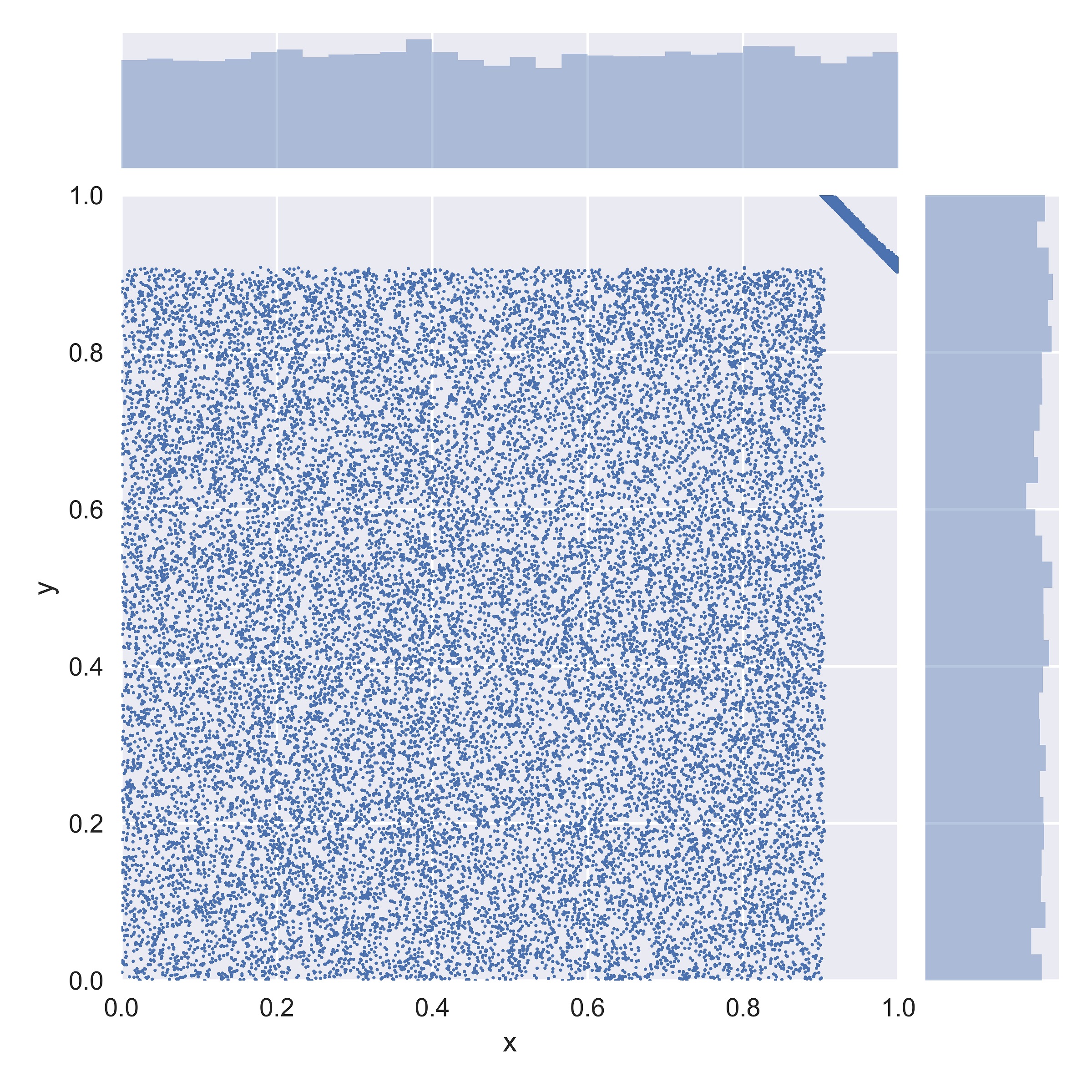}}
		\caption{Bounds on the distribution of the sum of dependent variables: 
			Sampled points from the numerically optimal measure $\nu^*$ and 
			the corresponding empirical 
			marginal distributions.}
		\label{FigureRearr}
	\end{center}
\end{figure}

For simplicity, we consider the case $d = 2$, and $\mu_i = 
\mathcal{U}([0,1])$. Let $s = 1.9$. Every optimal measure 
$\nu \in \mathcal{Q}$ gives mass $1/10$ uniformly to the line section $\{(x, 1.9-x) : 0.9 
\leq x 
\leq 1 \}$, while the rest of the mass is irrelevant as long as the marginal 
condition is satisfied. This leads to an optimal value $\phi(f) = 0.1$. 
For the natural 
choice of reference measure $\theta = \mathcal{U}([0,1]^2)$, every optimal 
measure is singular with respect to $\theta$, and thus one can expect high 
errors by penalization. An implementation with this reference measure and $L^2$ 
penalization with $\gamma = 320$ leads to $\hat{\phi}_{\theta,\gamma}^m(f) 
\approx  
0.0881$.\footnote{The network was trained for 20000 iterations with batch size 
	1024.}

\textsl{Updating the reference measure:} To obtain a more accurate value, we 
make use of Equation \eqref{dualoptimizer} in Theorem \ref{thm:reg}. 
Recall that an optimizer $\hat{\nu} \in \mathcal{Q}$ of 
$\phi_{\theta,\gamma}(f)$ is given by $\frac{d\hat{\nu}}{d\theta} = 
\beta'(f-\hat{h})$, where $\hat{h}$ is the dual optimizer of 
$\phi_{\theta,\gamma}(f)$. Taking $\hat{\nu}$ as a reference 
measure 
instead of $\theta$ can only reduce the error by penalization, since 
$\phi_{\hat{\nu}, 
	\gamma}(f) \geq \phi_{\theta, \gamma}(f)$ holds by convexity of $\beta^\ast_{\gamma}$.
Implementing the problem with $\hat{\nu}$ as a reference measure has to be done 
approximately, since the true optimizer $\hat{h}$ is unknown and replaced by 
the numerical optimal solution. We denote the numerically obtained optimal 
measure by $\nu^*$. 

To implement $\phi_{\nu^*,\gamma}^m(f)$ requires sampling points from $\nu^*$. 
This is non-trivial since $\nu^*$ is only given by $\frac{d\nu^*}{d\theta}$. We 
implemented this by an acceptance-rejection method as described in 
\cite{flury1990acceptance}. This is very slow, as the number of rejections increases with the maximum of the Radon-Nikodym derivative. Sampling efficiently in such a situation is difficult, 
see e.g.~\cite{horger2018deep} for an overview of existing methods and a 
proposed new one.

Figure \ref{FigureRearr} illustrates the optimal 
measure 
$\nu^*$. The measure $\nu^*$ looks comparable to an optimal solution of 
$\phi(f)$, while simultaneously being driven towards the reference measure 
$\theta$. 
One obtains $\hat{\phi}_{\nu^*,\gamma}^m(f) \approx 0.0982$, which is close to 
the true optimal value $0.1$.

In the following, we briefly discuss the rearrangement
algorithm \cite{embrechts2013model, puccetti2012computation} which is tailored 
to this type of problem. In contrast to the presented 
approach, which relies on sampling from the involved marginal distributions, 
the rearrangement algorithm mainly relies on the (inverse of) the cumulative 
distribution function of the marginals. The rearrangement 
algorithm achieves similar or even better accuracy 
in higher dimensional settings and with different marginals (e.g.~Pareto 
marginals). The case of higher dimensions scales well with the approach taken here. However, the base time 
in low dimensions is higher than that of the rearrangement algorithm, and 
further heavy tailed 
marginals like the Pareto distribution can lead to less 
accuracy.\footnote{We refer to 
\url{https://github.com/stephaneckstein/transport-and-related}, where different 
specifications (including Pareto marginals) of the problem in this section are 
implemented.}

\subsection{Generative Adversarial Networks (GANs)}
The objective in GANs is to create new sample 
points from a measure $\mu$, of which only an empirical distribution 
$\tilde{\mu}$ is known (see e.g.~\cite{arjovsky2017wasserstein, goodfellow2014generative}). Usually, the measure $\mu$ might refer to the uniform distribution over some very large set of images. The set $\tilde{\mu}$ is then just the uniform distribution over a small subset of these images. The goal is to sample new images that are not already present in the given subset, but that might plausibly have been samples from the measure $\mu$.

To proceed, we first take some latent probability measure $\tau$. The goal is to 
obtain a function $G$ such that $\mu$ and the push-forward measure $\tau \circ G^{-1}$ are close (in a sense to be specified), and thus the pseudo samples for $\mu$ can be obtained by 
sampling from $\tau$ and applying $G$. 
To find such a function $G$ out of a class of functions $\mathcal{G}$, one can only use $\tilde\mu$ instead of $\mu$ (thus $\mu$ does not enter the formal problem statement). The closeness of $\tilde\mu$ and 
$\tau \circ G^{-1}$ is measured by different distances in GANs. In the 
Wasserstein GAN, the first Wasserstein distance $W_1(\cdot,\cdot)$ (see e.g~\cite{villani2008optimal})
is used, and the
objective is
\[
\argmin_{G \in \mathcal{G}} W_1(\tau \circ G^{-1}, \tilde\mu).
\]

The above can be generalized to arbitrary transport distances instead of 
the first Wasserstein distance. To put this into our setting, let $\mathcal{G}$ a set of functions that map into $\mathcal{X}_1$. Let $\mathcal{X} = \mathcal{X}_1 
\times \mathcal{X}_1$ and for $G\in\mathcal{G}$ define $\mathcal{Q}_G := \{\nu 
\in \mathcal{P}(\mathcal{X}) : \nu_1 = \tau \circ G^{-1},~\nu_2 = \tilde\mu\}$. For a cost 
function $c$, arbitrary transport type GANs can be expressed via
\[
\argmin_{G \in \mathcal{G}} - \phi_G(-c)
= \argmin_{G \in \mathcal{G}} \inf_{\nu \in \mathcal{Q}_G} \int c \,d\nu
\]
If $c$ is a metric, the above corresponds to the Wasserstein GAN.

Since it is difficult to objectively evaluate GAN setups, we omit numerical 
results in this section. The interested reader can see the code on GitHub, 
where the toy problems appearing in \cite{gulrajani2017improved} are implemented
for different functions $c$.\footnote{For the implementation, we 
merely adjusted code from 
\url{https://github.com/igul222/improved_wgan_training} to our method. See again
\url{https://github.com/stephaneckstein/transport-and-related}} 

\section{Proofs}\label{sec5}
\subsection{Proof of Theorem \ref{thm:reg}}
{\rule{0mm}{1mm}\\[-6ex]\rule{0mm}{1mm}}

1) We first show \eqref{eq:dualreg} by verifying that
$\phi_{\theta,\gamma}$ is real-valued and continuous from above on $C_b(\mathcal{X})$. To that end, as shown in Appendix \ref{sec:DS} one has $\phi^\ast(\mu)=\sup_{h\in\mathcal{H}}(\int h\,d\mu-\int h\,d\mu_0)$ for all $\mu\in\mathcal{P}(\mathcal{X})$, so that   
\begin{equation}\label{rep:Q}
\mu\in\mathcal{Q}\quad\mbox{if and only if}\quad \int h\,d\mu=\int h\,d\mu_0\mbox{ for all }h\in\mathcal{H}.
\end{equation} 
Since $\beta_\gamma(x)\ge xy-\frac{1}{\gamma}\beta^\ast(y)$ for all $x\in\mathbb{R}$ and $y\in\mathbb{R}_+$, it follows that
\[
\int\beta_\gamma(f-h)\,d\theta \ge \int f-h\,d\pi -\frac{1}{\gamma}\int\beta^\ast\Big(\frac{d\pi}{d\theta}\Big)\,d\theta
\]
so that 
\[
\phi_{\theta,\gamma}(f)=\inf_{h\in\mathcal{H}}\Big\{\int h\,d\pi+\int \beta_\gamma(f-h)\,d\theta\Big\} \ge \int f\,d\pi -\frac{1}{\gamma}\int\beta^\ast\Big(\frac{d\pi}{d\theta}\Big)\,d\theta>-\infty
\]
for all $f\in C_b(\mathcal{X})$. This shows that $\phi_{\theta,\gamma}$ is real-valued on $C_b(\mathcal{X})$. Further, let
$(f^k)$ be a sequence in $C_b(\mathcal{X})$ such that $ f^k\downarrow f$. For every $h\in\mathcal{H}$ the monotone convergence theorem implies $\int \beta_\gamma(f^k-h)\,d\theta\to \int \beta_\gamma(f-h)\,d\theta$,
so that $\phi_{\theta,\gamma}(f^k)\downarrow \phi_{\theta,\gamma}(f)$. Hence, it 
follows from the nonlinear Daniell-Stone theorem (see Proposition \ref{prop:nonlinearDS}) that
\[\phi_{\theta,\gamma}(f)=\max_{\mu\in \mathcal{P}(\mathcal{X})}\Big\{\int 
f\,d\mu-\phi^{\ast}_{\theta,\gamma}(\mu) \Big\}\quad\mbox{for all }f\in C_b(\mathcal{X}),\]
where the convex conjugate is given by 
\[\phi^{\ast}_{\theta,\gamma}(\mu)=\phi^\ast(\mu)+\psi^\ast_{\theta,\gamma}(\mu)=\begin{cases} \frac{1}{\gamma}\int 
\beta^\ast\big(\frac{d\mu}{d\theta}\big)\,d\theta & \mbox{if }\mu\in\mathcal{Q}\mbox{ and  }\mu\ll\theta \\ 
\infty &\mbox{else}.\end{cases}\] 
Indeed, the convex conjugate of the convolution $\inf_{f\in C_b(\mathcal{X})}\{\phi(f)+\psi_{\theta,\gamma}(\cdot-f)\}$ is given as the sum of the convex conjugates $\phi^\ast$ and $\psi_{\theta,\gamma}^\ast$. 
By \eqref{rep:Q} one has $\phi^\ast(\mu)=0$ if $\mu\in\mathcal{Q}$ and $\phi^\ast(\mu)=+\infty$ otherwise. Moreover, 
\[
\psi^\ast_{\theta,\gamma}(\mu)=\sup_{f\in C_b(\mathcal{X})}\Big\{\int f\,d\mu-\int \beta_\gamma(f)\,d\theta\Big\}=\sup_{f\in C_b(\mathcal{X})}\Big\{\int f\frac{d\mu}{d\theta}- \beta_\gamma(f)\,d\theta\Big\}=\int \beta_\gamma^\ast\Big(\frac{d\mu}{d\theta}\Big)\,d\theta
\] 
if $\mu\ll\theta$ and $\psi^\ast_{\theta,\gamma}(\mu)=+\infty$ otherwise.

2) We next show \eqref{eq:boundmain}. On the one hand, one has
\[
\phi_{\theta,\gamma}(f)=\inf_{h\in\mathcal{H}}\Big\{\int h\,d\mu_0+\int \beta_\gamma(f-h)\,d\theta\Big\}\le \inf_{\substack{h\in\mathcal{H}:\\h\ge f}}\int h\,d\mu_0+\beta_\gamma(0)=\phi(f)+\frac{\beta(0)}{\gamma}.
\]
On the other hand, for every $\varepsilon$-optimizer $\mu_\varepsilon\in\mathcal{Q}$ of \eqref{superhedgingdualneu} such that $\mu_\epsilon\ll\theta$ one has
\begin{align*}
\phi(f)&\le \int f\,d\mu_\varepsilon+\varepsilon
\le \int f\,d\mu_\varepsilon-\phi^\ast_{\theta,\gamma}(\mu_\varepsilon)+\phi^\ast_{\theta,\gamma}(\mu_\varepsilon)+\varepsilon\le\phi_{\theta,\gamma}(f)+\frac{1}{\gamma}\int\beta^\ast\Big(\frac{d\mu_\varepsilon}{d\theta}\Big)\,d\theta+\varepsilon
\end{align*} 
with the convention $-\infty+\infty=+\infty$.

3) Let $\hat h\in \mathcal{H}$ be a minimizer of \eqref{defphithetagamma}, 
i.e.~$\phi_{\mu,\gamma}(f)=\int \hat h\,d\mu_0+
\int \beta_\gamma(f-\hat h)\, d\theta$. Defining $h_\lambda:=\hat h+\lambda h$ for an arbitrary 
$h\in\mathcal{H}$, the first order condition
\[\frac{d}{d\lambda}\Big|_{\lambda=0}\Big(\int h_\lambda\,d\mu_0+\int \beta_\gamma(f- h_\lambda)\, 
d\theta \Big)=0\]
implies
\[\int h\,d\mu_0-\int \beta^\prime_\gamma(f-\hat h) h\,d\theta  = 0.\]
This shows that the probability measure $\hat\mu$ with Radon-Nikod\'ym derivative $\frac{d\hat\mu}{d\theta}:=\beta_{\gamma}^\prime(f-\hat h)$ satisfies $\int h\,d\mu_0=\int h\,d\hat\mu$ for all $h\in\mathcal{H}$, which in view of \eqref{rep:Q} satisfies $\hat\mu\in\mathcal{Q}$. Integrating the identity 
$\beta_\gamma(x)=x\beta_\gamma^\prime(x)-\beta_\gamma^\ast(\beta_\gamma^\prime(x))$
with $x=f-\hat h$ w.r.t.~$\theta$, one obtains
\[
\int \beta_\gamma(f-\hat h)\,d\theta = \int f-\hat h\,d\hat\mu-\int\beta^\ast_\gamma\Big(\frac{d\hat\mu}{d\theta}\Big)\,d\theta
\]
which shows that
\[
\phi_{\theta,\gamma}(f)=\int\hat h\,d\mu_0+\int\beta_\gamma(f-\hat h)\,d\theta=\int f\,d\hat\mu-\int\beta^\ast_\gamma\Big(\frac{d\hat\mu}{d\theta}\Big)\,d\theta.
\]
As a consequence, $\hat\mu\in\mathcal{Q}$ is a maximizer of \eqref{eq:dualreg}.
\qed

\subsection{Proof of Proposition \ref{prop:approx}}
Fix $f\in C_b(\mathcal{X})$.
That $\lim_{m\to\infty}\phi^m(f)=\phi^\infty(f)\ge \phi(f)$ follows 
from the definition of $\mathcal{H}^\infty$. Moreover, for every 
$\varepsilon>0$ the Condition (D) guarantees 
$h\in\mathcal{H}^\infty$ and 
$K\subseteq\mathcal{X}$ such that $1_{K^c}\le h$ and $\int 
h\,d\mu_0\le\varepsilon$. Hence, $\phi^\infty(1_{K^c})\leq \int h 
\,d\mu_0 \le\varepsilon$ and Dini's lemma implies that $\phi^\infty$ is continuous from above on $C_b(\mathcal{X})$. By Proposition \ref{prop:nonlinearDS} it follows that
\[\phi^\infty(f)=\max_{\mu\in \mathcal{P}(\mathcal{X})} \Big\{\int 
f\,d\mu-\phi^{\infty\ast}(\mu)\Big\}.\] Similar to \eqref{conj:phi:kappa} its convex conjugate is given by
\[\phi^{\infty\ast}(\mu)=\sup_{h\in\mathcal{H}^\infty}\Big(\int h\,d\mu- \int h\,d\mu_0\Big)\le 
\sup_{h\in\mathcal{H}}\Big(\int h\,d\mu- \int h\,d\mu_0\Big)=\phi^\ast(\mu).\]
It remains to show that for $h\in\mathcal{H}$ and $\mu\in \mathcal{P}(\mathcal{X})$ 
with $\int h\,d\mu- \int h\,d\mu_0>0$ there exists 
$h^\prime\in\mathcal{H}^\infty$ such that $\int h^\prime\,d\mu- \int h^\prime\,d\mu_0>0$. 
But this follows directly from the first part of Condition 
(D) for the probability measure $\frac{1}{2}\mu+\frac{1}{2}\mu_0$. Indeed,
there 
exists a sequence $(h^n)$ in $\mathcal{H}^\infty$ such that $h^n\to h$ in 
$L^1(\mu)$ and in $L^1(\mu_0)$, which shows that $\int h^n\,d\mu- \int h^n\,d\mu_0>0$ for $n$ large enough. 
\qed

\subsection{Proof of Proposition \ref{prop:unif}}
Observe that	
\[
\phi^m(f) + \beta_{\gamma}(0) \geq \inf_{\substack{h \in \mathcal{H}^m:\\h\geq f}} \Big\{\int h\,d\mu_0 + \int \beta_\gamma(f-h) \,d\theta\Big\} \geq \phi^m_{\theta, \gamma}(f) \geq \phi_{\theta, \gamma}(f)
\]
where the first inequality uses that $\beta_\gamma$ is increasing, the second inequality just drops the constraint $h\geq f$, and the third inequality follows from $\mathcal{H}^m \subseteq \mathcal{H}$. 

Fix $\varepsilon>0$. By Condition (D) and Theorem \ref{thm:reg} there exist  $m_0 \in \mathbb{N}$ and $\gamma_0>0$ such that  
\[
\phi^m(f)\le \phi(f)+\varepsilon\quad\mbox{and}\quad \phi(f)\le\phi_{\theta,\gamma}(f)+\varepsilon
\]
for all $m\ge m_0$ and $\gamma\ge\gamma_0$. This shows that 
\[
\phi(f) + \varepsilon +\frac{\beta(0)}{\gamma}\geq \phi^m(f) +\beta_\gamma(0) \geq \phi_{\theta,\gamma}^m(f) \geq \phi_{\theta,\gamma}(f) \geq \phi(f) -\varepsilon
\]
for all $m\ge m_0$ and $\gamma\ge\gamma_0$, which shows that $\phi^m_{\theta, \gamma}(f) \rightarrow \phi(f)$ whenever $\min\{m, \gamma\} \rightarrow \infty$. 
\qed

\subsection{Proof of Lemma \ref{lem:CondD}}
\label{proof:CondD}
\begin{itemize}
	\item[(a)] From Hornik \cite{hornik1991approximation} it follows that 
	$\mathfrak{N}_{l_j, d_j}$ is dense in $C_b(\mathbb{R}^{d_j})$ with respect 
	to 
	$L^1(\nu)$ for
	every $\nu \in \mathcal{P}(\mathbb{R}^{d_j})$ and all $j = 1, ..., J$. By 
	the 
	triangle inequality and boundedness of $e_j$, 
	the first part of 
	Condition $(D)$ follows.
	\item[(b)] If $\mathcal{X}$ is compact, the condition is trivially 
	satisfied. Hence assume that $\mathcal{X} = \mathbb{R}^d = 
	\mathbb{R}^{d_1} \times ... \times \mathbb{R}^{d_{J_0}}$ and $\pi_j = {\rm pr}_j,~ 
	e_j 
	= 1$ for $j=1,..., J_0 \leq J$, where ${\rm pr}_j$ is the projection 
	from 
	$\mathbb{R}^d$ to the $j$-th marginal component in 
	$\mathbb{R}^{d_j}$.
	
	Let $\varepsilon > 0$ and $\nu \in \mathcal{P}(\mathcal{X})$. We first fix $j$, denote 
	by $\nu^{(j)} := \nu \circ {\rm pr}_j^{-1}$ and show 
	that there exists a 
	$h^{\prime\prime}_j \in 
	\mathfrak{N}_{l_j, d_j}$ such that $1_{K_j^c} \leq h^{\prime\prime}_j$ and 
	$\int_{\mathbb{R}^{d_j}} 
	h^{\prime\prime}_j\,
	d\nu^{(j)} \leq 2\varepsilon$ for some 
	compact subset $K_j$ of $\mathbb{R}^{d_j}$. 
	Without loss of generality, assume that $l_j = 1$. This can always be done 
	since 
	the function $h^{\prime\prime}_j$ will be compact-valued and hence for 
	multiple layers, 
	the remaining layers beyond the first can simply approximate the identity 
	function in the supremum norm.\footnote{More precisely, the function 
		$h^{\prime\prime}_j$ as given by one layer would be the input in the 
		first 
		component 
		for the remaining layers, and the remaining layers approximate the 
		continuous function $[-z,z]^m \ni x \mapsto x_1$ in the supremum norm, 
		which is possible as shown in \cite{hornik1991approximation}.} Fix 
	$K_j=[-c,+c]^{d_j}$ such that 
	$\nu^{(j)}(K_j^c)\le\varepsilon/(4d_j)$.
	By assumption on $\varphi$, for each $i\in\{1,\dots,d_j\}$ there exist 
	$\underline a_i,\underline b_i,\overline a_i,\overline b_i\in\mathbb{R}$ 
	such 
	that 
	\[\varphi(\underline a_i x_i+\underline b_i)+\varphi(\overline a_i 
	x_i+\overline b_i)\begin{cases} \le \varepsilon/(2d_j) & \mbox{for 
	}x_i\in[-c,c] 
	\\ \ge 1-\varepsilon & \mbox{for }x_i\not\in[-c-1,c+1]. \end{cases}\]
	Then
	\[h^{\prime\prime}_{j}:=\sum_{i=1}^{d_j} \varphi(\underline a_i 
	x_i+\underline 
	b_i)+\varphi(\overline a_i x_i+\overline b_i)+\varepsilon 
	\in\mathfrak{N}_{1,d_j, 2d_j}\subset \mathfrak{N}_{1, d_j} \]
	satisfies $1_{\tilde K_j^c} \leq h^{\prime\prime}_{j} $ for the compact 
	$\tilde 
	K_j:=[-c-1,c+1]^{d_j}$, as well as
	\[\int_{\mathbb{R}^{d_j}} h^{\prime\prime}_{j}\,d\nu^{(j)}\le 
	\int_{\mathbb{R}^{d_j}}  
	\frac{\varepsilon}{2} 1_{K_j} + 
	2d 
	1_{K_j^c}\,d\nu^{(j)} + \varepsilon \le 
	\frac{\varepsilon}{2}+2d\nu^{(j)}(K_j^c) +
	\varepsilon \le 2\varepsilon.\]
	
	Now, define $h^{\prime\prime} := \sum_{j=1}^d h^{\prime\prime}_j \circ {\rm 
		pr}_j \in \mathcal{H}^{\infty}$ and $K := \prod_{j=1}^{J_0} \tilde{K}_j 
	\subset 
	\mathcal{X}$, which is compact. Then one immediately gets $1_{K^c} \leq 
	h^{\prime\prime}$ and $\int h^{\prime\prime} d\nu \leq 2 J_0 \varepsilon$.
\end{itemize}
\subsection{Proof of Proposition \ref{nouniform}}
1) For one fix network $\mathfrak{N}_{l_j, d_j, m}$, 
the mapping $\xi \mapsto N_{l_j, d_j, m}(\xi)$ is pointwise 
continuous, 
i.e.~it holds for $\xi_n \rightarrow \xi$ that $N_{l_j, d_j, 
	m}(\xi_n)(x) \rightarrow N_{l_j, d_j, m}(\xi)(x)$ 
for all $x\in \mathbb{R}^{d_j}$, 
since $\varphi$ is continuous. Further, since we assume that $\varphi$ is 
bounded, the functions $N_{l_j, d_j, m}(\xi_n)$ are uniformly bounded 
and 
hence by 
dominated convergence one 
obtains $N_{l_j, d_j, m}(\xi_n) \rightarrow N_{l_j, d_j, m}(\xi)$ in 
$L^1(\nu)$ for all $\nu \in 
\mathcal{P}(\mathbb{R}^{d_j})$. By the triangle inequality, this continuity 
transfers to the mapping $(\xi_1, ..., \xi_J) \mapsto \sum_{j=1}^J 
e_j N_{l_j, d_j, m}(\xi_j) \circ \pi_j$. Hence, we can write 
\[
\mathcal{H}^m = \{ \eta(A) + a : A \in \mathcal{A}^m, a\in \mathbb{R}\}
\]
where $A \mapsto \eta(A)$ is continuous in $L^1(\mu_0)$ and $L^1(\theta)$, and $\mathcal{A}^m$ is compact.

2) 
For every $\varepsilon>0$ there exists 
$\eta(A)+a\in\mathcal{H}^m$ with 
$\eta(A)+a\ge f$ such that 
\[\phi^m(f)+\varepsilon\ge \int\eta(A)+a\,d\mu_0=\lim_{\gamma\to\infty} 
\Big\{\int\eta(A)\,d\mu_0+ a + \int \beta_\gamma(f-\eta(A)-a)\,d\theta 
\Big\}\ge\limsup_{\gamma\to\infty}\phi^m_{\theta,\gamma}(f)\]
since $0\le \int \beta_\gamma(f-\eta(A)-a)\,d\theta\le 
\beta_\gamma(0)=\frac{1}{\gamma}\beta(0)$.

On the other hand, let $(\gamma_n)$ be a sequence in $\mathbb{R}_+$ with 
$\gamma_n\to\infty$. Our goal is to show that
$\phi^m(f)\le \liminf_{n\to\infty} \phi^m_{\theta,\gamma_n}(f)$. 
We 
assume 
that $ \liminf_{n\to\infty} \phi^m_{\theta,\gamma_n}(f)<\infty$ otherwise 
there is 
nothing to prove. For every $n\in\mathbb{N}$ 
there exist  $A^n\in\mathcal{A}^m$ and $a^n\in\mathbb{R}$ such that 
\begin{equation}\label{bound1}
\phi^m_{\theta,\gamma_n}(f) +\frac{1}{n} \ge \int\eta(A^n)\,d\mu_0+a^n + \int 
\beta_{\gamma_n}\big(f-\eta(A^n)-a^n\big)\,d\theta \ge \int\eta(A^n)\,d\mu_0 +
\int  f-\eta(A^n)+ c\,d\theta 
\end{equation}
since $\beta_{\gamma_n}(x)\ge x+c$ for all $n\in\mathbb{N}$ for some 
constant 
$c\in\mathbb{R}$.
In particular, $\liminf_{n\to\infty} \phi^m_{\theta,\gamma_n}(f)$ is 
real-valued 
as $A\mapsto \int  f-\eta(A)+ c\,d\theta$ is a continuous function on the 
compact 
$\mathcal{A}^m$. By passing to a subsequence we may assume that 
$\lim_{n\to\infty} \phi^m_{\theta,\gamma_n}(f)=\liminf_{n\to\infty} 
\phi^m_{\theta,\gamma_n}(f)$, and  $A^n\to A\in \mathcal{A}^m$ such that 
$\eta(A^n)\to \eta(A)$ in $L^1(\theta)$ and $\theta$-a.s.~as well as 
$\int\eta(A^n)\,d\mu_0\to \int\eta(A)\,d\mu_0$. We next show that $(a^n)$ is bounded. 
Suppose 
by way of contradiction that $a^n\to-\infty$. Since  $\lim_{x\to\infty} 
\beta(x)/x =\infty$ and $f - \eta(A^n)$ is uniformly bounded by compactness 
of $\mathcal{A}^m$, it follows 
that
\[\int\eta(A^n)\,d\mu_0+ a^n + 
\beta_{\gamma_n}\big(f-\eta(A^n)-a^n\big)\to+\infty\] 
Moreover, in view of \eqref{bound1} the sequence \[\big(\int\eta(A^n)\,d\mu_0+ a^n 
+ 
\beta_{\gamma_n}\big(f-\eta(A^n)-a^n\big)\big)^-\] is uniformly integrable 
in 
$L^1(\theta)$. Hence, it follows from Fatou's lemma that
\begin{align*}
+\infty&=\int \liminf_{n\to\infty} \big\{\int\eta(A^n)\,d\mu_0+ a^n 
+\beta_{\gamma_n}\big(f-\eta(A^n)-a^n\big)\big\}\,d\theta\\&\le 
\liminf_{n\to\infty}\Big\{ \int\eta(A^n)\,d\mu_0+ a^n  + \int  
\beta_{\gamma_n}\big(f-\eta(A^n)-a^n\big)\,d\theta\Big\}\\
&\le \liminf_{n\rightarrow \infty} \phi^m_{\theta,\gamma_n}(f) < \infty
\end{align*}
which is the desired contradiction. This shows that $(a^n)$ is bounded and 
by 
passing to a subsequence 
$a^n\to a\in\mathbb{R}$. Finally it follows from Fatou's lemma that
\begin{align*}
\liminf_{n\to\infty} \phi^m_{\theta,\gamma_n}(f)&= 
\liminf_{n\to\infty}\Big\{ 
\int\eta(A^n)\,d\mu_0+ a^n + \int  
\beta_{\gamma_n}\big(f-\eta(A^n)-a^n\big)\,d\theta\Big\} \\
&\ge \int\eta(A)\,d\mu_0+ a+ \int  \beta_{\infty}\big(f-\eta(A)-a\big)\,d\theta
\\&= \int\eta(A)+ a\,d\mu_0\\& \ge \phi^m(f)
\end{align*}
where $\beta_\infty(x)=0$ if $x\le 0$ and $\beta_\infty(x)=\infty$ if $x> 
0$.
The second inequality follows because $\eta(A)+a\ge f$ $\theta$-a.s.~as a 
consequence of the first inequality, $f,\eta(A)\in C_b(\mathcal{X})$ and $\theta$ is strictly positive.
\qed

\begin{appendix}
	\section{Nonlinear version of the Daniell-Stone theorem}\label{sec:DS}
	Let $\mathcal{X}$ be a Polish space. Given a measurable function $\kappa:\mathcal{X}\to [1,\infty)$, we denote by $C_\kappa(\mathcal{X})$ the Stone vector lattice of all continuous functions $f\colon\mathcal{X}\to\mathbb{R}$  such that $|f|/\kappa$ is bounded. For instance, if $\kappa$ is bounded one has  $C_\kappa(\mathcal{X})=C_b(\mathcal{X})$, or if $\kappa(x)=1+|x|$ on $\mathcal{X}=\mathbb{R}^d$ the space $C_\kappa(\mathbb{R}^d)$ contains all continuous functions $f\colon\mathbb{R}^d\to\mathbb{R}$ of linear growth. Further, let $ca^+_\kappa(\mathcal{X})$ be the set of all Borel measures $\mu$ on $\mathcal{X}$ which satisfy $\int \kappa\,d\mu<\infty$. The following nonlinear version of the Daniell-Stone theorem follows directly from  Proposition 1.1 in \cite{cheridito2015representation}.
	\begin{proposition}\label{prop:nonlinearDS}
		Let  $\phi\colon C_\kappa(\mathcal{X})\to\mathbb{R}$ be an increasing\footnote{$\phi(f)\ge\phi(g)$ whenever $f\ge g$.} convex 
		functional which is continuous from above, i.e.~$\phi(f^n)\downarrow 0$ 
		for every sequence $(f^n)$ such that $f^n\downarrow 0$. 
		Then, it has the dual representation
		\begin{equation}\label{rep:dual:kappa}
		\phi(f)=\max_{\mu\in ca^+_\kappa(\mathcal{X})} \Big\{\int f\,d\mu-\phi^\ast(\mu)\Big\}\quad\mbox{for 
			all }f\in C_\kappa(\mathcal{X}),
		\end{equation}
		where the convex conjugate 
		$\phi^\ast\colon ca^+_\kappa(\mathcal{X})\to\mathbb{R}\cup\{+\infty\}$ is given by 
		$\phi^\ast(\mu)=\sup_{f\in C_\kappa(\mathcal{X})} \{\int f\,d\mu-\phi(f)\}$. 	
	\end{proposition}
	Continuity from above is strongly related to the concept of  tightness, which in the context of risk measures was introduced 
	by F\"ollmer and Schied, see \cite{follmer2011stochastic}.
	Typical examples include transport type problems where tightness is imposed 
	by marginal constraints, see e.g.~Bartl et al.~\cite{bartl2017duality}. For 
	extensions of the representation \eqref{rep:dual:kappa} to upper 
	semicontinuous functions and related
	pricing-hedging dualities we refer to Cheridito et 
	al.~\cite{cheridito2017duality}.
	
	As an application we consider the superhedging functional
	\[
	\phi(f):=\inf\Big\{\int h\,d\mu_0: h\ge f\mbox{ for some 
	}h\in\mathcal{H}\Big\}
	\]
	on $C_\kappa(\mathcal{X})$, where $\mu_0\in ca^+_\kappa(\mathcal{X})$ is a probability measure and $\mathcal{H}\subseteq C_\kappa(\mathcal{X})$ is a convex cone such that $\kappa\in\mathcal{H}$. Straightforward inspection shows that $\phi$ is a real-valued increasing convex functional on $C_\kappa(\mathcal{X})$. Further, if $\phi$ is continuous from above by Proposition \ref{prop:nonlinearDS} it has the dual representation \eqref{rep:dual:kappa}.
	Its convex conjugate is given by
	\begin{align}
	\phi^\ast(\mu)&=\sup_{f\in C_\kappa(\mathcal{X})}\Big\{\int f\,d\mu-\inf_{\substack{h\in\mathcal{H}:\\h\ge f}}\int h\,d\mu_0\Big\}\nonumber\\
	&=\sup_{h\in \mathcal{H}}\sup_{\substack{f\in C_\kappa(\mathcal{X}):\\h\ge f}}\Big\{\int f\,d\mu-\int h\,d\mu_0\Big\}\nonumber\\
	&=\sup_{h\in \mathcal{H}}\Big\{\int h\,d\mu-\int h\,d\mu_0\Big\}.\label{conj:phi:kappa}
	\end{align}
	Since $\mathcal{H}$ is a convex cone which contains the constants it follows that $\phi^\ast(\mu)=0$ whenever $\mu\in ca^+_\kappa(\mathcal{X})$ is a probability measure such that $\int h\,d\mu=\int h\,d\mu_0$ for all $h\in\mathcal{H}$, and $\phi^\ast(\mu)=+\infty$ else. In particular, in case that $C_\kappa(\mathcal{X})=C_b(\mathcal{X})$ we conclude the dual representation \eqref{superhedgingdualneu}. 
\end{appendix}

\bibliographystyle{abbrv}

\end{document}